\PassOptionsToPackage{dvipsnames}{xcolor} 
\documentclass[11pt,letter]{article}
\usepackage{amssymb,amsmath,amsthm,multirow,color,xcolor,cancel,bbm,tikz,mathrsfs,hyperref,authblk,stmaryrd}
\usepackage[textheight=24cm,textwidth=16cm]{geometry}
\tikzstyle{main node}=[draw,circle,inner sep=1,outer sep=2,thick,minimum size=12pt]
\usetikzlibrary{decorations.pathreplacing}

\hypersetup{
    colorlinks=true,
    linkcolor=blue,
    citecolor=cyan,
    filecolor=magenta,      
    urlcolor=green,
    pdftitle={Overleaf Example},
    pdfpagemode=FullScreknownen,
    }

\newtheorem{proposition}{Proposition}
\newtheorem{lemma}{Lemma}
\newtheorem{theorem}{Theorem}
\newtheorem{corollary}{Corollary}

\newtheorem{remark}{Remark}
\newtheorem{example}{Example}

\newcommand{\Q}[1]{\llbracket #1\rrbracket}
\newcommand{\EM}[1]{{\it\textcolor{Maroon}{#1}}}
\newcommand{\BF}[1]{{\boldmath{\bf #1}\unboldmath}}
\newcommand{\B}{\{0,1\}}
\newcommand{\TWO}{\mathbf{2}}
\newcommand{\G}{\mathbb{G}}
\newcommand{\cst}{\mathrm{cst}}

\newcommand{\lcm}{\mathrm{lcm}}
\def\Im{\mathrm{Im}}

\title{Interaction graphs of isomorphic automata networks II:\\ universal dynamics}

\author{
Florian Bridoux\footnote{\scriptsize Universit\'e Côte d’Azur, CNRS, I3S, Sophia Antipolis, France. ({\tt  bridoux@i3s.unice.fr})},
Aymeric Picard Marchetto\footnote{\scriptsize Universit\'e Côte d’Azur, CNRS, I3S, Sophia Antipolis, France. ({\tt  picard@i3s.unice.fr})},
Adrien Richard\footnote{\scriptsize Universit\'e Côte d’Azur, CNRS, I3S, Sophia Antipolis, France. ({\tt  adrien.richard@cnrs.fr})},
}

\begin{document}

\maketitle

\begin{abstract}
An automata network with $n$ components over a finite alphabet $Q$ of size $q$ is a discrete dynamical system described by the successive iterations of a  function $f:Q^n\to Q^n$. In most applications, the main parameter is the  interaction graph of $f$: the digraph with vertex set $[n]$ that contains an  arc from $j$ to $i$ if $f_i$ depends on input $j$. What can be said on the  set $\G(f)$ of the interaction graphs of the automata networks isomorphic to  $f$? It seems that this simple question has never been studied. In a previous paper, we prove that the complete digraph $K_n$, with $n^2$ arcs, is universal in that $K_n\in \G(f)$ whenever $f$ is not constant nor the identity (and $n\geq 5$). In this paper, taking the opposite direction, we prove that there exist universal automata networks $f$, in that $\G(f)$ contains all the digraphs on $[n]$, excepted the empty one. Actually, we prove that the presence of only three specific digraphs in $\G(f)$ implies the universality of $f$, and we prove that this forces the alphabet size $q$ to have at least $n$ prime factors (with multiplicity). However, we prove that for any fixed $q\geq 3$, there exists almost universal functions, that is, functions $f:Q^n\to Q^n$ such that the probability that a random digraph belongs to $\G(f)$ tends to $1$ as $n\to\infty$. We do not know if this holds in the binary case $q=2$, providing only partial results. 
\end{abstract}

\section{Introduction}

An \EM{automata network} with $n$ components over the finite alphabet of size $q\geq 2$ is a function 
\[
f:Q^n\to Q^n,\quad x=(x_1,\dots,x_n)\mapsto f(x)=(f_1(x),\dots,f_n(x))
\]
where $Q=\{0,1,\dots,q-1\}$. The components $f_i:Q^n\to Q$ are usually called the \EM{local transition functions} of the network, while $f$ is referred as the \EM{global transition function} (but here we identify this function as the network). Automata networks are also called \EM{finite dynamical systems}, and in the binary case, when $Q=\{0,1\}$, the term \EM{Boolean networks} is applied. We denote by \EM{$F(n,q)$} the set of automata networks with $n$ components over the alphabet of size $q$.

\medskip
The dynamics described by $f$ is explicitly represented by the digraph \EM{$\Gamma(f)$} with vertex set $Q^n$ and an arc from $x$ to $f(x)$ for every $x\in Q^n$. Two automata networks $f,h:Q^n\to Q^n$ are \EM{isomorphic}, $f\sim h$ in notation, if $\Gamma(f)$ and $\Gamma(h)$ are isomorphic in the usual sense. An equivalence class of $\sim$ then corresponds to an unlabeled digraph with $q^n$ vertices in which each vertex has out-degree exactly one. 

\medskip
Automata networks have many applications. In particular, they are omnipresent in the modeling of neural and gene networks (see \cite{N15} for a review). The ``network'' terminology comes from the fact that the \EM{interaction graph} of $f$ is often considered as the main parameter: it is the digraph \EM{$G(f)$} with vertex set $[n]=\{1,\dots,n\}$ and such that, for all $i,j\in [n]$, there is an arc from $j$ to $i$ if $f_i$ depends on input $j$, that is, if there exist $x,y\in Q^n$ which only differ in $x_j\neq y_j$ such that $f_i(x)\neq f_i(y)$. 

\medskip
In many applications, as modeling of gene networks, the interaction graph is often well approximated while the actual dynamics is not \cite{TK01,N15}. One is thus faced with the following question: {\em what can be said on the dynamics described by $f$ from $G(f)$ only?} There are many results in this direction; see \cite{G20} for a review. In most cases, the studied dynamical properties are invariant by isomorphism: number of fixed points, number of images, number of periodic configurations, number and size of limit cycles, transient length and so on. However, the interaction graph is {\em not} invariant by isomorphism: even if $f$ and $h$ are isomorphic, their interaction graphs can be very different. This variation can give some {\em limit} to the central question stated above and we think it deserves some study.

\medskip
Hence, we propose a study of this variation, considering the set \EM{$\G(f)$} of interaction graphs of automata networks isomorphic to $f$:
\[
\G(f)=\{G(h)\mid h\in F(n,q),~h\sim f\}.
\]
If $G\in\G(f)$ we say that $f$ \EM{can be produced} by $G$. For instance, if $f$ is constant, then $\G(f)$ contains a single
digraph, the digraph on $[n]$ without arcs, and if $f$ is the identity, then $\G(f)$ also contains a single digraph, the digraph on $[n]$ with $n$ loops (cycles of length one) and no other arcs. So the identity and constant automata networks can be produced by a single digraph.

\medskip
In a previous paper \cite{bppmr2023}, we proved, with Perrot, the following. First, for $n\geq 5$, $f$ can always be produced by the complete digraph on $[n]$ (with $n^2$ arcs), excepted if $f$ is the identity or a constant; and the complete digraph on $[n]$ is the unique digraph with this property. Second, for $q$ fixed and $n$ large enough, $\G(f)$ contains a digraph distinct from the complete digraphs. Together, this shows that, for $q$ fixed and $n$ large enough, the identity and constant automata networks are the only networks which can be produced by a single digraph. This gives some insight on automata networks $f$ minimizing $|\G(f)|$. 

\medskip
In this paper, we take the converse direction: we study automata networks $f$ for which $|\G(f)|$ is large. In this case, $f$ is produced by many digraphs, and thus the dynamics described by $f$ ``weakly'' depends on the interaction graph of $f$. Obviously, $|\G(f)|$ is at most $2^{n^2}$, the number of digraphs on $[n]$, but we can be slightly more precise:
\[
|\G(f)|\leq 2^{n^2}-1. 
\]
Indeed, if $\G(f)=2^{n^2}$, then $\G(f)$ contains the digraph on $[n]$ without arcs, which implies that $f$ is constant, but then $\G(f)$ only contains the empty digraph, a contradiction. 

\medskip
We will prove that this bound is tight, and we will call the automata networks $f$ reaching this bound \EM{universal}, in that they can be produced by all the digraphs excepted the empty one. That universal automata networks exist is not obvious, and we actually prove that $f$ is universal whenever $\G(f)$ contains three very specific digraphs. 

\begin{theorem}\label{thm:universal_1_intro}
An automata network $f\in F(n,q)$ is universal whenever $\G(f)$ contains an acyclic digraph of height~$1$, a digraph with one loop and no other arcs, and a digraph with $n$ loops and no other arcs. 
\end{theorem}

Using this, we then characterize the couples $(n,q)$ such that $F(n,q)$ contains a universal automata network; this shows that $q$ must have at least $n$ primes factors (with multiplicity). 

\begin{theorem}\label{thm:universal_2_intro}
For all $n,q\geq 2$, $F(n,q)$ contains a universal automata network if and only if $q$ can be expressed as the product of $n$ integers $q=q_1q_2\dots q_n$ with $\prod_{i=1}^n(q_i-1)\geq 2^n$. 
\end{theorem}

As a corollary, we obtain that $F(n,3^n)$ contains a universal automata network, and that if $F(n,q)$ contains a universal function, then $q$ must be exponential with $n$, namely $q\geq 3^n$. So for $n$ fixed, $q=3^n$ is the smallest alphabet size allowing a universal automata network to exist. 

\medskip
But using an alphabet of size exponential in $n$ takes us away from the spirit of automata networks: network components are supposed ``simple'', working on a small alphabet, the ``complexity'' coming from the large number $n$ of components and their interconnections. To deal with this more usual situation, we fix the alphabet size $q$, and study the behavior of   
\[
\gamma_q(n)=\max_{f\in F(n,q)} |\G(f)|,
\]
which is the maximum number of digraphs on $[n]$ that can produce the same $f\in F(n,q)$. If $\gamma_q(n)/2^{n^2}\to 1$ as $n\to\infty$, it means that, for some $f\in F(n,q)$, the probability that a random digraph on $[n]$ can produce $f$ tends to $1$ as $n\to\infty$. We then say that $f$ can be produced by almost all digraphs. We prove that this happens for every non-binary alphabet $q$. Actually, since almost all digraphs are Hamiltonian \cite{W73}, this follows from the following result. 

\begin{theorem}\label{thm:Hamiltonian1}
For all $n\geq 1$ and $q\geq 3$, there exists $f\in F(n,q)$ such that $\G(f)$ contains all the Hamiltonian digraphs on $[n]$. 
\end{theorem}

We provide another family of automata networks $f$ produced by almost all digraphs, using a lemma that should be of independent interest. Given a digraph $G$ on $[n]$ and $q\geq 2$, let us say that $G$ is \EM{$q$-permutable} if $F(n,q)$ contains a permutation whose interaction graph is $G$, and let us say that $G$ is \EM{coverable} if its vertices can be spanned by vertex disjoint cycles. Gadouleau \cite{G18rank} proved that every $q$-permutable digraph is $G$ coverable, and that the converse holds for $q\geq 3$ but fails for $q=2$. Our lemma says that, for $q=2$, the converse fails not ``too much''. 

\begin{lemma}\label{lem:2-permutable_intro}
Almost all coverable digraphs are $2$-permutable. 
\end{lemma}

In the binary case $q=2$, we have almost no result concerning $\gamma_2(n)$ and analyzing the asymptotic behaviors of $\gamma_2(n)/2^{n^2}$ seems difficult. We were only able to prove, using again the above lemma, the following bounds: for $\epsilon>0$ and $n$ large enough, 
\begin{equation}\label{eq:gamma2}
(1-\epsilon)2^{(n-1)^2}\leq \gamma_2(n)\leq 2^{n^2}-2^{\frac{n^2}{2}}.
\end{equation}

\medskip
The difficulty in analyzing $\gamma_2(n)$ leads us to study a modified notion of universality. Given an integer $k\geq 1$, we say that $f\in F(n,2)$ is \EM{$k$-induced universal} if, for every $k$-vertex digraph $H$, some digraph in $\G(f)$ has an induced subgraph isomorphic to $H$. It is then interesting to consider the smallest integer $n(k)$ such that $F(n(k),2)$ contains a $k$-induced universal automata network. We show that $n(k)$ can be estimated rather accurately (logarithms are in base~$2$):
\begin{equation}\label{eq:Nk}
k+\log k-\log(\lceil \log k\rceil+1)-1\leq n(k)\leq k+\lceil \log k\rceil+1.
\end{equation}

\medskip
The paper is organized as follows. Definitions and basic results are given in Section \ref{sec:def}. Universality is studied in Section \ref{sec:universal_functions}; Theorems~\ref{thm:universal_1_intro} and \ref{thm:universal_2_intro} are proven there. The asymptotic behaviors of $\gamma_q(n)$ in the non-binary case $q\geq 3$ is studied in Section \ref{sec:almost_universal}; Theorem~\ref{thm:Hamiltonian1} and Lemma~\ref{lem:2-permutable_intro} are proven there. The binary case is studied in Section~\ref{sec:binary}; inequalities \eqref{eq:gamma2} and \eqref{eq:Nk} are proven there. Finally, concluding remarks are given in Section \ref{sec:conclu}.  

\section{Basic definitions}\label{sec:def}

In this section, we collect basic definitions, including those from the introduction. 

\medskip
Throughout the paper, $n$ and $q$ are always integers, with $n\geq 1$ and $q\geq 2$. We set $\EM{[n]}=\{1,\dots,n\}$ and $\EM{\Q{q}} = \{0,1,\dots,q-1\}$. Given a finite set $I$, we denote by \EM{$\Q{q}^I$} the set of functions $x$ from $I$ to $\Q{q}$. Such a function $x$ is called a \EM{configuration} on the set of \EM{components} $I$. Given $i\in I$, we write \EM{$x_i$} instead of $x(i)$; this is the \EM{state} of component $i$ in configuration $x$. Given $J\subseteq I$, we denote by \EM{$x_J$} the restriction of $x$ on $J$. We write \EM{$\Q{q}^{n}$} instead of $\Q{q}^{[n]}$, and we then represent $x\in\Q{q}^n$ as a row vector $x=(x_1,\dots,x_n)$. We denote by \EM{$0^n$} the configuration with $n$ zeroes, and \EM{$1^n$} the configuration with $n$ ones. We denote by \EM{$e_i$} the configuration such that $(e_i)_i=1$ and $(e_i)_j=0$ for all $j\neq i$. Given two configurations $x,y$ on the same set of components $I$, and $i\in I$, we say that $x$ and $y$ \EM{only differ in $i$} if $x_i\neq y_i$ and $x_j=y_j$ for all $j\in I\setminus\{i\}$. Logarithms are always in base~$2$, and we sometime use \EM{$\oplus$} for the addition modulo $2$.

\medskip
We denote by \EM{$F(n,q)$} the set of functions from $\Q{q}^n$ to itself. These functions are usually called automata networks with $n$ components over the alphabet of size $q$, but we simply call them functions in the following. Given $f\in F(n,q)$ and $i\in [n]$ we denote by \EM{$f_i$} the function from $\Q{q}^n$ to $\Q{q}$ such that $f_i(x)=f(x)_i$ for all $x\in\Q{q}^n$. We denote by \EM{$\Im(f)$} the set of images of $f$, and $|\Im(f)|$ is the \EM{rank} of $f$. A \EM{fixed point} of $f$ is a configuration $x\in\Q{q}^n$ such that $f(x)=x$. Given $k\geq 1$, we denote by \EM{$f^k$} the $k$ fold composition of $f$. We say that $f$ is \EM{$k$-nilpotent} if $f^k=\cst$ but $f^{k-1}\neq \cst$; such a function has a unique fixed point. Two functions $f,h\in F(n,q)$ are \EM{isomorphic} (or \EM{conjugate}) is there exists a permutation $\pi$ of $\Q{q}^n$ such that $\pi\circ f=h\circ \pi$; then $\pi$ is an isomorphism from $f$ to $h$ (and $\pi^{-1}$ is an isomorphism from $h$ to $f$). An equivalence class of $\sim$ then corresponds to an unlabeled digraph with $q^n$ vertices in which each vertex has out-degree exactly one.  

\medskip
A digraph $G$ on a finite set $V$ is a digraph with vertex set $V$. A \EM{source} of $G$ is a vertex of in-degree $0$, and a \EM{sink} is a vertex of out-degree $0$. The \EM{height} of an acyclic digraph $G$ is the number of arcs in a longest path of $G$ from a source to a sink. We denote by \EM{$C_n$} the directed cycle of length $n$ whose vertices are $1,2,\dots,n$ in order. Given two digraphs $G,H$, we write $G\,\EM{\subseteq}\, H$ to mean that $G$ is a subgraph of $H$.

\medskip
The \EM{interaction graph} of $f\in F(n,q)$, denoted \EM{$G(f)$}, is the digraph on $[n]$ such that, for all $i,j\in [n]$, there is an arc from $i$ to $j$ whenever there exists two configurations $x,y\in\Q{q}^{n}$ which only differ in $x_j\neq y_j$ such that $f_i(x)\neq f_i(y)$. Given a digraph $G$ on $[n]$, we denote by \EM{$F(G,q)$} the set of functions $f\in F(n,q)$ with $G(f)=G$. We denote by \EM{$\G(f)$} the set of interaction graphs of functions isomorphic to $f$:
\[
\G(f)=\{G(h)\mid h\in F(n,q),~h\sim f\}.
\]
If $G\in\G(f)$ we say that $f$ \EM{can be produced} by $G$. Two digraphs $G,H$ on $[n]$ are isomorphic if there exists a permutation $\pi$ of $[n]$ such that, for all $i,j\in [n]$, $G$ has an arc from $i$ to $j$ if and only if $H$ has an arc from $\pi(i)$ to $\pi(j)$; we say that $\pi$ is an isomorphism from $G$ to $H$. We stress that $\G(f)$ is \EM{closed under isomorphism}: if $G\in \G(f)$ and $H$ is a digraph on $[n]$ isomorphic to $G$, then $H\in\G(f)$. Indeed, let $\pi$ be an isomorphism from $G$ to $H$, and let $\tilde \pi$ be the permutation of $\Q{q}^n$ defined by $\tilde\pi(x)=(x_{\pi(1)},\dots,x_{\pi(n)})$ for all $x\in\Q{q}^n$. If $G\in\G(f)$, there exists $g\in F(G,q)$ isomorphic to $f$, and $H$ is the interaction graph of $h=\tilde\pi^{-1}\circ g\circ \tilde \pi$, thus $H\in \G(f)$ since $h\sim g$ and $g\sim h$. 
  
\section{Universal functions}\label{sec:universal_functions}

We say that a function $f\in F(n,q)$ is \EM{universal} if it can be produced by any digraph on $[n]$ excepted the empty one, that is, if $|\G(f)|$ reaches the trivial bounds $2^{n^2}-1$. An equivalent formulation is that $f\in F(n,q)$ is universal if, for any digraph $G$ on $[n]$ with at least one arc, $F(G,q)$ contains a function isomorphic to $f$. 

\medskip
In this section, we show that such functions exist, providing a rather simple and explicit characterization. We further give necessary and sufficient conditions on $n$ and $q$ for $F(n,q)$ to contain a universal function.   

\subsection{Universal functions are $2$-nilpotent}\label{sec:2-nilpotent}

A fundamental result concerning automata networks is Robert's theorem. 

\begin{theorem}[Robert \cite{R80,R86}]\label{thm:robert}
If $f\in F(n,q)$ and $G(f)$ is acyclic of height $k$, then $f^{k+1}=\cst$.
\end{theorem}

As an immediate consequence, we get the following. 

\begin{lemma}\label{lem:2-nilpotent}
If $f\in F(n,q)$ is universal, then $f$ is $2$-nilpotent.
\end{lemma}

\begin{proof}
Let $G$ be an acyclic digraph on $[n]$ of height $1$. If $f\in F(n,q)$ is universal, there exists $h\in F(G,q)$ isomorphic to $f$. By Robert's theorem, $h^2=\cst$, and since $G$ has at least one arc, $h\neq\cst$. Thus $h$ is $2$-nilpotent  and since $h\sim f$, $f$ is also $2$-nilpotent. 
\end{proof}

We deduce that there are no universal functions in $F(n,2)$ because $F(G,2)$ does not contain $2$-nilpotent functions for some $G$ with at least one arc: for instance, by the following lemma, if $G$ is a disjoint union of cycles, then $F(G,2)$ only contains permutations, and thus it does not contain a $2$-nilpotent function. 

\begin{lemma}\label{lem:disjoint_union_of_cycles}
If $f\in F(n,2)$ and $G(f)$ is a disjoint union of cycles, then $f$ is a permutation. 
\end{lemma}

\begin{proof}
That $G(f)$ is a disjoint union of cycles means that there exists a permutation $\pi$ of $[n]$ such that, for all $i\in [n]$, $\pi(i)$ is the unique in-neighbor of $i$ in $G$. Then, for all $i\in [n]$ there exists $a_i\in\B$ such that $f_i(x)=x_{\pi(i)}\oplus a_i$ for all $x\in\B^n$. Let $x,y\in\B^n$ distinct. Then there exists $i\in [n]$ such that $x_{\pi(i)}\neq y_{\pi(i)}$ and we obtain $f_i(x)=x_{\pi(i)}\oplus a_i\neq y_{\pi(i)}\oplus a_i=f_i(y)$. Thus $f(x)\neq f(y)$, that is, $f$ is a permutation. 
\end{proof}

\begin{remark}
For alphabets of size $q\geq 3$ the situation changes drastically: $F(G,q)$ contains a $2$-nilpotent function whenever $G$ has at least one arc. Indeed, let $f\in F(G,q)$ defined as follows: 
\[
f_i(x)=
\left\{
\begin{array}{l}
1\textrm{ if $x_j\geq 2$ for some in-neighbor $j$ of $i$ in $G$,}\\[1mm]
0\textrm{ otherwise.}
\end{array}
\right.
\]
Then $\Im(f)\subseteq \B^n$ and $f(x)=0^n$ for all $x\in\B^n$, thus $f$ is $2$-nilpotent (see \cite{GR16} for other results concerning the existence of nilpotent functions). But the isomorphic class of $f$ depends on $G$, thus $f$ is not  universal: for instance, if $G$ is a disjoint union of cycles, then $f$ has $2^n$ images, while if $G$ is the complete digraph (with $n^2$ arcs) then $f$ has only $2$ images, $0^n$ and $1^n$. 
\end{remark}

\begin{figure}
\[
\begin{array}{ccc}
\begin{array}{c}
\begin{tikzpicture}
\useasboundingbox (-1.9,-0.5) rectangle (2,2.4);
\def\s{0.10}
\def\v{0.2}
\def\is{-1}
\def\os{-1}
\node[circle,inner sep=\is,outer sep=\os] (n1) at (-1.5,1) {\tiny$\bullet$};
\node[circle,inner sep=\is,outer sep=\os] (n2) at (-0.5,1) {\tiny$\bullet$};
\node[circle,inner sep=\is,outer sep=\os] (n3) at (0.5,1) {\tiny$\bullet$};
\node[circle,inner sep=\is,outer sep=\os] (n0) at (0,0) {\tiny$\bullet$};
\node[circle,inner sep=\is,outer sep=\os] (n01) at ({1.5-0.2},1) {\tiny$\bullet$};
\node[circle,inner sep=\is,outer sep=\os] (n02) at ({1.5},1) {\tiny$\bullet$};
\node[circle,inner sep=\is,outer sep=\os] (n03) at ({1.5+0.2},1) {\tiny$\bullet$};
\node[circle,inner sep=\is,outer sep=\os] (n11) at ({-1.5-0.3},2) {\tiny$\bullet$};
\node[circle,inner sep=\is,outer sep=\os] (n12) at ({-1.5-0.1},2) {\tiny$\bullet$};
\node[circle,inner sep=\is,outer sep=\os] (n13) at ({-1.5+0.1},2) {\tiny$\bullet$};
\node[circle,inner sep=\is,outer sep=\os] (n14) at ({-1.5+0.3},2) {\tiny$\bullet$};
\node[circle,inner sep=\is,outer sep=\os] (n21) at ({-0.5-0.2},2) {\tiny$\bullet$};
\node[circle,inner sep=\is,outer sep=\os] (n22) at ({-0.5},2) {\tiny$\bullet$};
\node[circle,inner sep=\is,outer sep=\os] (n23) at ({-0.5+0.2},2) {\tiny$\bullet$};
\node[circle,inner sep=\is,outer sep=\os] (n31) at ({0.5-0.1},2) {\tiny$\bullet$};
\node[circle,inner sep=\is,outer sep=\os] (n32) at ({0.5+0.1},2) {\tiny$\bullet$};
\node[circle,inner sep=0.7,outer sep=0.7] (n0b) at (0,0) {};
\draw[->] (n0b.-112) .. controls ({0-0.2},-0.3) and ({0+0.2},-0.3) .. (n0b.-68);
\draw 
(n1) edge (n0)
(n2) edge (n0)
(n3) edge (n0)
(n01) edge (n0)
(n02) edge (n0)
(n03) edge (n0)
(n11) edge (n1)
(n12) edge (n1)
(n13) edge (n1)
(n14) edge (n1)
(n21) edge (n2)
(n22) edge (n2)
(n23) edge (n2)
(n31) edge (n3)
(n32) edge (n3)
;
\end{tikzpicture}
\end{array}
&\qquad\qquad&
\begin{array}{c}
\begin{tikzpicture}
\useasboundingbox (-1.9,-0.5) rectangle (2,2.4);
\def\s{0.10}
\def\v{0.2}
\def\is{-1}
\def\os{-1}
\draw[fill=gray!20] (-1.5,1) -- ({-2+2*\s},2) -- ({-1-2*\s},2) -- (-1.5,1);
\draw[fill=gray!20] (-0.5,1) -- ({-1+2*\s},2) -- ({0-2*\s},2) -- (-0.5,1);
\draw[fill=gray!20] (0.5,1) -- ({0+2*\s},2) -- ({1-2*\s},2) -- (0.5,1);
\draw[fill=gray!20] (0,0) -- ({1+\s},1) -- ({2-\s},1) -- (0,0);
\node[circle,inner sep=\is,outer sep=\os] (n1) at (-1.5,1) {\tiny$\bullet$};
\node[circle,inner sep=\is,outer sep=\os] (n2) at (-0.5,1) {\tiny$\bullet$};
\node[circle,inner sep=\is,outer sep=\os] (n3) at (0.5,1) {\tiny$\bullet$};
\node[circle,inner sep=\is,outer sep=\os] (n0) at (0,0) {\tiny$\bullet$};
\node[circle,inner sep=0.7,outer sep=0.7] (n0b) at (0,0) {};
\draw[->] (n0b.-112) .. controls ({0-0.2},-0.3) and ({0+0.2},-0.3) .. (n0b.-68);
\draw 
(n1) edge (n0)
(n2) edge (n0)
(n3) edge (n0)
;
\node  at (-1.5,2+2*\s) {\scriptsize$4$};
\node  at (-0.5,2+2*\s) {\scriptsize$3$};
\node  at (+0.5,2+2*\s) {\scriptsize$2$};
\node  at (1.5,1+2*\s) {\scriptsize$3$};
\end{tikzpicture}
\end{array}\\
\text{(a)} && \text{(b)}
\end{array}
\]
{\caption{\label{fig:2nil} (a) An unlabelled $2$-nilpotent function $f\in F(4,2)$; edges are directed toward the loop. (b) Schematic representation (used several times in the following).}}
\end{figure}
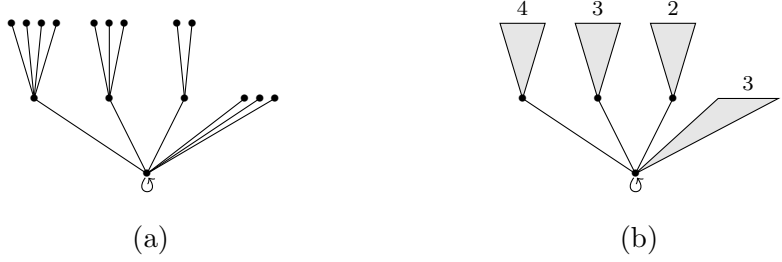

$2$-nilpotent functions have a very simple structure: up to isomorphism, they are completely described by the number of images, the number of pre-images of its fixed point, and the number of pre-images of the other images (see Figure~\ref{fig:2nil}). Testing isomorphism is then easy. 

\begin{lemma}\label{lem:2-nil_iso}
Let $f,h\in F(n,q)$ be $2$-nilpotent. Then $f$ and $h$ are isomorphic if and only if there exists a bijection $\pi:\Im(f)\to\Im(h)$ which maps the fixed point of $f$ to the fixed point of $h$ and which preserves the number of pre-images, that is, $|f^{-1}(x)|=|h^{-1}(\pi(x))|$ for all $x\in\Im(f)$. 
\end{lemma}

\begin{proof}
Let $\pi$ be as in the statement. Suppose, without loss, that $0^n$ is the fixed point of $f$; hence $\pi(0^n)$ is the fixed point of $h$. Since $f,h$ are $2$-nilpotent, we have $\Im(f)\subseteq f^{-1}(0^n)$ and $\Im(h)\subseteq h^{-1}(\pi(0^n))$, and since $|f^{-1}(x)|=|h^{-1}(\pi(x))|$ for all $x\in\Im(f)$, we deduce that there is a permutation $\tilde \pi$ of $\Q{q}^n$ extending $\pi$ such that $\tilde \pi(f^{-1}(x))=h^{-1}(\pi(x))$ for all $x\in\Im(f)$. Then $\tilde \pi$ is an isomorphism between $f$ and $h$: for any $x\in \Q{q}^n$, we have $\tilde\pi(x)\in h^{-1}(\pi(f(x)))$ and thus $h(\tilde\pi(x))=\pi(f(x))=\tilde \pi(f(x))$. This proves one direction and the other is obvious. 
\end{proof}

\subsection{Universal functions exist}

In this section we prove the existence of universal functions. More precisely, we prove the existence of universal \EM{regular} functions, where $f\in F(n,q)$ is \EM{regular} if each image has the same number of pre-images, that is, there is $d$ such that $|f^{-1}(x)|=d$ for all $x\in\Im(f)$; $d$ is then the \EM{degree} of $f$. Our interest for this case is that universal regular functions are more easy to define and analyze than universal functions in general. Their study also allows us to introduce techniques that will be reused in the general case. 

\medskip
Note that, by Lemma~\ref{lem:2-nil_iso}, two regular $2$-nilpotent functions $f,h\in F(n,q)$ are isomorphic if and only if they have the same rank (since the preservation of the number of pre-images is given by the regularity).

\medskip
Our first result is a sufficient condition on $(n,q,r)$ for $F(n,q)$ to contain a universal regular function of rank $r$. In particular, it shows that if $q$ has at least $2n$ prime factors (with multiplicity), then $F(n,q)$ contains a universal regular function. 

\begin{lemma}\label{lem:regular_universal}
Let $n\geq 1$ and let $r\geq 2$ be an integer with at least $n$ prime factors (with multiplicity). Let $q\geq r^2$ be a multiple of $r$. Then $F(n,q)$ contains a universal regular function of rank $r$. 
\end{lemma}

\begin{proof}
Since $r$ has at least $n$ prime factors, there are $r_1,\dots,r_n\geq 2$ such that $r=r_1r_2\dots r_n$. Since $r\mid q$ and $q\geq r^2$, there is a set $Q\subseteq\mathbb{N}^{2n}$ of size $q$ such that 
\begin{equation}\label{eq:Q}
\Q{r_1}\times \cdots \times\Q{r_n}\times \Q{r_1}\times \dots \times\Q{r_n}\subseteq Q\subseteq \mathbb{N}^n\times \Q{r_1}\times \cdots \times\Q{r_n}.
\end{equation}
We then identify the alphabet $\Q{q}$ with $Q$. So given $f\in F(n,q)$ and $x\in Q^n$, the state $x_i$ of component $i$ at $x$ is $x_i=(x_{i,1},\dots,x_{i,2n})\in Q$ and the state of component $i$ at $f(x)$ is $f_i(x)=(f_{i,1}(x),\dots,f_{i,2n}(x))\in Q$. 

\medskip
Let $G$ be any digraph on $[n]$ with at least one arc. We will prove that $F(G,q)$ contains a $2$-nilpotent regular function of rank $r$. This proves the lemma since, by Lemma~\ref{lem:2-nil_iso}, any two $2$-nilpotent regular functions with the same rank are isomorphic. 

\medskip
Let $N(i)$ be the in-neighbors of $i$ in $G$, and let $J$ be the set of vertices of out-degree at least one in $G$ (which is non-empty since $G$ has at least one arc). Let $\phi$ be a function from $[n]$ to $J$ such that $\phi(j)=j$ for all $j \in J$. Consider the function $f\in F(n,q)$ defined as follows: for all $x\in Q^n$, $i\in [n]$ and $j\in [2n]$, 
\begin{equation}\label{eq:def_f}
f_{i,j}(x)=
\left\{
\begin{array}{ll}
x_{\phi(j),n+j} & \textrm{if }j\in [n]\textrm{ and }\phi(j) \in N(i),\\
0 & \textrm{otherwise.} 
\end{array}
\right.
\end{equation}
Note that $f_i(x)\in \Q{r_1}\times \cdots \times\Q{r_n}\times \{0\}\times \cdots \times \{0\}\in Q$ and thus $f$ is indeed a function from $Q^n$ to itself. See Example~\ref{ex:Q} for an illustration.

\medskip
We will prove that $G(f)=G$ and then that $f$ is $2$-nilpotent and regular with rank $r$. 

\medskip
Suppose that $G$ has an arc from $j$ to $i$. Let $x,y \in Q^n$ which only differ in $x_{j,n+j}\neq y_{j,n+j}$.
Since $j \in N(i)$ we have $j \in J$, and thus $\phi(j)=j$.  Hence 
\[
f_{i,j}(x) = x_{\phi(j),n+j} = x_{j,n+j} \neq y_{j,n+j} =  y_{\phi(j),n+j} = f_{ij}(y).
\]
Therefore $f_{i}(x) \neq f_i(y)$ and we deduce that $G(f)$ has an arc from $j$ to $i$. 
Conversely suppose that $G(f)$ has an arc from $j$ to $i$. Then there exists $x,y \in Q^n$, which only differ in $x_j\neq y_j$, such that $f_i(x)\neq f_i(y)$. So $f_{i,k}(x)\neq f_{i,k}(y)$ for some $k\in [n]$ with $\phi(k)\in N(i)$ and we obtain $x_{\phi(k),n+k}\neq y_{\phi(k),n+k}$. Since $x_{\phi(k)}\neq y_{\phi(k)}$ we have $\phi(k)=j$ and we deduce that $G$ has an arc from $j$ to $i$. This proves that $G(f)=G$.


\medskip
We now prove that $f$ is $2$-nilpotent. Let $X$ be the set of $x\in Q^n$ such that: $x_{i,j}=0$ if $j\not\in [n]$ or $\phi(j) \notin N(i)$; and $x_{i,j}\in \Q{r_j}$ otherwise. The images of $f$ are all in $X$. Furthermore, for all $x\in X$ and $j\in [n]$ we have $x_{\phi(j),n+j}=0$, and thus $f(x)=(0^{2n})^n$. Hence $f^2=\cst=(0^{2n})^n$ and since $f\neq\cst$ this proves that $f$ is $2$-nilpotent. 

\medskip
We finally prove that $f$ is regular with rank $r$. First, remark that the value of $f(x)$ only depends on the values of $(x_{\phi(1),n+1}, x_{\phi(2),n+2}, \dots, x_{\phi(n),2n}) \in \Q{r_1}\times \cdots \times\Q{r_n}$. So the number of pre-images of $f(x)$ is a multiple of $q^n/r$, and $f$ has at most $r$ images. Now, consider two configurations $x,y \in Q^n$ such that  $x_{\phi(j),n+j}\neq y_{\phi(j),n+j}$ for some $j\in [n]$. Since $\phi(j) \in J$, there exist $i \in [n]$ such that $\phi(j) \in N(i)$ and we obtain $f_{i,j}(x) = x_{\phi(j),n+j} \neq y_{\phi(j),n+j} = f_{i,j}(y)$.
As a result, $f$ has at least $r$ images. Hence, $f$ has rank $r$, and this forces each image to have exactly $q^n/r$ pre-images.   
\end{proof}

\begin{example} \label{ex:Q}
Let $n = 3$ and $q=3^22^3$. Let $r_1=r_2=r_3=2$ and $r=r_1 r_2 r_3=2^3$. Then $q$ is a multiple of $r$ and $q=3^22^3=(2^3+1)2^3\geq r^2$. Hence, by Lemma~\ref{lem:regular_universal}, $F(n,q)$ contains a universal function of rank $r$. Let us illustrate the proof in this particular case. We can identify $\Q{q}$ with 
\[ 
Q = ( ( \Q{2}\times \Q{2}\times \Q{2}) \cup \{ (3,0,0) \} ) \times \Q{2}\times \Q{2}\times \Q{2},
\]
which satisfies $\Q{2}^6\subseteq Q\subseteq \mathbb{N}^3\times\Q{2}^3$. Now, let us illustrate that a regular function in $F(n,q)$ of rank $r$ can be produced by any non-empty digraph (this proves that a universal function in $F(n,q)$ exists by Lemma~\ref{lem:2-nil_iso}). Consider the following digraph $G$ on $[n]$: 
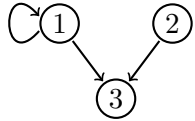
\[
\begin{tikzpicture}
\useasboundingbox (-0.7,-1.4) rectangle (2.2,0.4);
\node[outer sep=1,inner sep=2,circle,draw,thick] (1) at (0,0){$1$};
\node[outer sep=1,inner sep=2,circle,draw,thick] (2) at (1.5,0){$2$};
\node[outer sep=1,inner sep=2,circle,draw,thick] (3) at (0.75,-1){$3$};
\draw[->,thick] (1.{180+20}) .. controls (-0.8,-0.7) and (-0.8,+0.7) .. (1.{180-20});
\path[->,thick]
(2) edge (3)
(1) edge (3)
;
\end{tikzpicture}
\]
The set $J$ of vertices of out-degree at least one is $\{1,2\}$. So $\phi(1) = 1, \phi(2) = 2$ and we can choose $\phi(3)=1$. The corresponding function $f$, as defined in \eqref{eq:def_f}, is as follows:
\[
\begin{array}{|c|c|c|}
\multicolumn{1}{c}{x_1} & \multicolumn{1}{c}{x_2} & \multicolumn{1}{c}{x_3} \\
\hline
x_{1,1} & x_{2,1} & x_{3,1} \\
\hline
x_{1,2} & x_{2,2} & x_{3,2} \\
\hline
x_{1,3} & x_{2,3} & x_{3,3} \\
\hline
\textrm{\boldmath $x_{1,4}$\unboldmath} & x_{2,4} & x_{3,4} \\
\hline
x_{1,5} & \textrm{\boldmath $x_{2,5}$\unboldmath} & x_{3,5} \\
\hline
\textrm{\boldmath $x_{1,6}$\unboldmath} & x_{2,6} & x_{3,6} \\
\hline
\end{array}
\quad \xrightarrow{\text{f}} \quad
\begin{array}{|c|c|c|}
\multicolumn{1}{c}{f_1(x)} & \multicolumn{1}{c}{f_2(x)} & \multicolumn{1}{c}{f_3(x)} \\
\hline
\textrm{\boldmath $x_{1,4}$\unboldmath} & 0 & \textrm{\boldmath $x_{1,4}$\unboldmath} \\
\hline
0 & 0 & \textrm{\boldmath $x_{2,5}$\unboldmath} \\
\hline
\textrm{\boldmath $x_{1,6}$\unboldmath} & 0 & \textrm{\boldmath $x_{1,6}$\unboldmath} \\
\hline
0 & 0 & 0 \\
\hline
0 & 0 & 0 \\
\hline
0 & 0 & 0 \\
\hline
\end{array}
\quad \xrightarrow{\text{f}} \quad
\begin{array}{|c|c|c|}
\multicolumn{1}{c}{f^2_1(x)} & \multicolumn{1}{c}{f^2_2(x)} & \multicolumn{1}{c}{f^2_3(x)} \\
\hline
0 & 0 & 0 \\
\hline
0 & 0 & 0 \\
\hline
0 & 0 & 0 \\
\hline
0 & 0 & 0 \\
\hline
0 & 0 & 0 \\
\hline
0 & 0 & 0 \\
\hline
\end{array}
\]
Clearly, $G$ is the interaction graph of $f$, and $f$ is 2-nilpotent. Since $f(x)$ only depends on $(x_{1,4},x_{2,5},x_{1,6})\in\Q{2}^3$, the number of pre-images of $f(x)$ is a multiple of $q^n/r$, and $f$ has at most $r$ images. Furthermore, given $x,y\in Q^n$, if $(x_{1,4},x_{2,5},x_{1,6})\neq (y_{1,4},y_{2,5},y_{1,6})$ then $f(x)\neq f(y)$. Hence, $f$ has rank $r$, and this forces each image to have exactly $q^n/r$ pre-images.   
\end{example}

We now prove that the conditions of the previous lemma are necessary. More precisely, we prove that if $f\in F(n,q)$ is a $2$-nilpotent regular function of rank $r$ produced by only two specific digraphs, then the triple $(n,q,r)$ satisfies the conditions of the previous lemma. These two digraphs are: 
\begin{itemize}
\item
\EM{$L_{1,n}$}, the digraph on $[n]$ with a loop on vertex $1$ and no other arcs; and 
\item
\EM{$L_{n,n}$}, the digraph on $[n]$ with a loop on each vertex, and no other arcs.
\end{itemize}

\begin{lemma}\label{lem:rqn}
Let $f\in F(n,q)$ be a $2$-nilpotent regular function of rank $r$ with $L_{1,n},L_{n,n}\in\G(f)$. Then $r$ has at least $n$ prime factors (with multiplicity) and $r\mid q$ and $q\geq r^2$. 
\end{lemma}

\begin{proof}
First, since $L_{n,n}\in G(f)$, there is $g\in F(L_{n,n},q)$ isomorphic to $f$. Since $g_i$ only depends on component $i$ we have $\Im(g)=\Im(g_1)\times\cdots\times \Im(g_n)$. Setting $r_i=|\Im(g_i)|$ we have $r_i\geq 2$ since $g_i\neq\cst$, and we deduce that $r=|\Im(f)|=|\Im(g)|=r_1r_2\cdots r_n$ so $r$ has at least $n$ prime factors. 

\medskip
Second, since $L_{1,n}\in\G(f)$, there is $h\in F(L_{1,n},q)$ isomorphic to $f$. Thus $h$ is regular $2$-nilpotent of rank $r$ and we can suppose, without loss, that $0^n$ is the fixed point of $h$. Since $h_1$ only depends on component $1$, there is $h'_1\in F(1,q)$ such that $h_1(x)=h'_1(x_1)$ for every $x\in\Q{q}^n$. Since $h_i=\cst=0$ for every $i\neq 1$, we deduce that $\Im(h)=\Im(h'_1)\times\{0^{n-1}\}$, and thus $h'_1$ has rank $r$. Furthermore, for all $x\in\Q{q}^n$, we have $h^{-1}(x)={h'}^{-1}_1(x)_1\times \Q{q}\times\dots\times \Q{q}$. Since $h$ is regular, we deduce that $h'_1$ is regular, say with degree $d$. Then $q=rd$. Furthermore, since $h$ is $2$-nilpotent, $h'_1$ is $2$-nilpotent, and thus the $r$ images of $h'_1$ are pre-images of the fixed point of $h'_1$. Consequently, $d\geq r$ so that $q\geq r^2$. 
\end{proof}

Putting things together, we obtain the following characterizations. 

\begin{theorem}\label{thm:universal_regular_1}
Let $n,q\geq 2$ and let $f\in F(n,q)$ be regular. The following conditions are equivalent:
\begin{itemize}
\item
$f$ is universal.
\item
$f$ is $2$-nilpotent and $L_{1,n},L_{n,n}\in\G(f)$.
\end{itemize}

\end{theorem}

\begin{proof}
Suppose that $f$ satisfies the second point. By Lemma~\ref{lem:rqn}, the rank $r$ of $f$ has at least $n$ prime factors (with multiplicity) and $r\mid q$ and $q\geq r^2$. Consequently, by Lemma~\ref{lem:regular_universal}, $F(n,q)$ contains a universal regular function $h$ of rank $r$. By Lemma~\ref{lem:2-nil_iso}, $f$ and $h$ are isomorphic, and thus $f$ is universal. The other direction is obvious using Lemma~\ref{lem:2-nilpotent}. 
\end{proof}

\begin{theorem}\label{thm:universal_regular_2}
For $n,q\geq 2$, the following conditions are equivalent:
\begin{itemize}
\item
$F(n,q)$ contains a universal regular function of rank $r$.
\item
there is $r\in\mathbb{N}$ with at least $n$ prime factors (with multiplicity) such that $r\mid q$ and $q\geq r^2$. 
\end{itemize}
\end{theorem}

\begin{proof}
Suppose that $f\in F(n,q)$ is universal and regular. By Lemma~\ref{lem:2-nilpotent}, $f$ is $2$-nilpotent, and since $f$ is universal, we have $L_{1,n},L_{n,n}\in \G(f)$. Hence, by Lemma~\ref{lem:rqn}, the rank $r$ of $f$ has at least $n$ prime factors (with multiplicity) and $r\mid q$ and $q\geq r^2$. This proves one direction and the other is given by Lemma~\ref{lem:regular_universal}. 
\end{proof}

\begin{remark}
Theorem \ref{thm:universal_regular_2} shows that $F(n,4^n)$ contains a universal regular function and that $F(n,q)$ do not contain such a function for $q<4^n$, so $4^n$ is the minimal alphabet for the existence of a universal regular function with $n$ components. Furthermore, if $f\in F(n,4^n)$ is universal and regular, then, by Lemma \ref{lem:rqn}, its rank is necessarily $2^n$. We then deduce from Lemma~\ref{lem:2-nil_iso} that, up to isomorphism, $F(n,4^n)$ contains a unique universal regular function, which is thus the minimal universal regular function for $n$ components. 
\end{remark}

\begin{example}
For $n=2$, any universal regular function $f\in F(2,16)$ is $2$-nilpotent with $4$ images, each with $64$ pre-images. It has thus the following schematic representation (see Figure \ref{fig:2nil}):
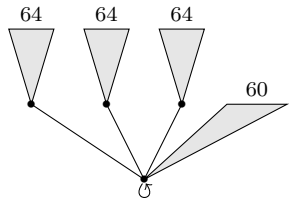
\[
\begin{tikzpicture}
\def\s{0.10}
\def\v{0.2}
\def\is{-1}
\def\os{-1}
\draw[fill=gray!20] (-1.5,1) -- ({-2+2*\s},2) -- ({-1-2*\s},2) -- (-1.5,1);
\draw[fill=gray!20] (-0.5,1) -- ({-1+2*\s},2) -- ({0-2*\s},2) -- (-0.5,1);
\draw[fill=gray!20] (0.5,1) -- ({0+2*\s},2) -- ({1-2*\s},2) -- (0.5,1);
\draw[fill=gray!20] (0,0) -- ({1+\s},1) -- ({2-\s},1) -- (0,0);
\node[circle,inner sep=\is,outer sep=\os] (n1) at (-1.5,1) {\tiny$\bullet$};
\node[circle,inner sep=\is,outer sep=\os] (n2) at (-0.5,1) {\tiny$\bullet$};
\node[circle,inner sep=\is,outer sep=\os] (n3) at (0.5,1) {\tiny$\bullet$};
\node[circle,inner sep=\is,outer sep=\os] (n0) at (0,0) {\tiny$\bullet$};
\node[circle,inner sep=0.7,outer sep=0.7] (n0b) at (0,0) {};
\draw[->] (n0b.-112) .. controls ({0-0.2},-0.3) and ({0+0.2},-0.3) .. (n0b.-68);
\draw 
(n1) edge (n0)
(n2) edge (n0)
(n3) edge (n0)
;
\node  at (-1.5,2+2*\s) {\scriptsize$64$};
\node  at (-0.5,2+2*\s) {\scriptsize$64$};
\node  at (+0.5,2+2*\s) {\scriptsize$64$};
\node  at (1.5,1+2*\s) {\scriptsize$60$};
\end{tikzpicture}
\] 
The universality can be easily checked with Theorem~\ref{thm:universal_regular_1}. For that, we have to show that $L_{1,2}$ and $L_{2,2}$ are in $\G(f)$. That $L_{2,2}\in\G(f)$ means that $f$  is isomorphic to the product of two $2$-nilpotent functions in $F(1,16)$, and indeed:
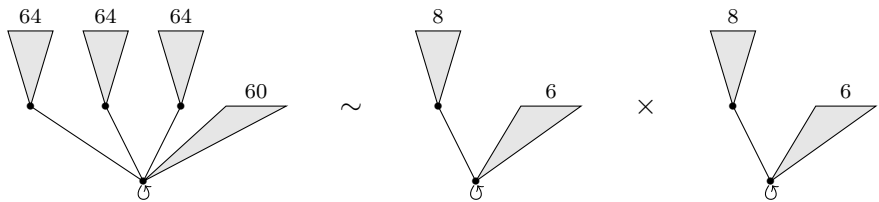
\[
\begin{array}{ccccc}
\begin{array}{c}
\begin{tikzpicture}
\def\s{0.10}
\def\v{0.2}
\def\is{-1}
\def\os{-1}
\draw[fill=gray!20] (-1.5,1) -- ({-2+2*\s},2) -- ({-1-2*\s},2) -- (-1.5,1);
\draw[fill=gray!20] (-0.5,1) -- ({-1+2*\s},2) -- ({0-2*\s},2) -- (-0.5,1);
\draw[fill=gray!20] (0.5,1) -- ({0+2*\s},2) -- ({1-2*\s},2) -- (0.5,1);
\draw[fill=gray!20] (0,0) -- ({1+\s},1) -- ({2-\s},1) -- (0,0);
\node[circle,inner sep=\is,outer sep=\os] (n1) at (-1.5,1) {\tiny$\bullet$};
\node[circle,inner sep=\is,outer sep=\os] (n2) at (-0.5,1) {\tiny$\bullet$};
\node[circle,inner sep=\is,outer sep=\os] (n3) at (0.5,1) {\tiny$\bullet$};
\node[circle,inner sep=\is,outer sep=\os] (n0) at (0,0) {\tiny$\bullet$};
\node[circle,inner sep=0.7,outer sep=0.7] (n0b) at (0,0) {};
\draw[->] (n0b.-112) .. controls ({0-0.2},-0.3) and ({0+0.2},-0.3) .. (n0b.-68);
\draw 
(n1) edge (n0)
(n2) edge (n0)
(n3) edge (n0)
;
\node  at (-1.5,2+2*\s) {\scriptsize$64$};
\node  at (-0.5,2+2*\s) {\scriptsize$64$};
\node  at (+0.5,2+2*\s) {\scriptsize$64$};
\node  at (1.5,1+2*\s) {\scriptsize$60$};
\end{tikzpicture}
\end{array}
&\begin{array}{c}\sim\end{array}&
%
\begin{array}{c}
\begin{tikzpicture}
\def\s{0.10}
\def\v{0.2}
\def\is{-1}
\def\os{-1}
\draw[fill=gray!20] (-0.5,1) -- ({-1+2*\s},2) -- ({0-2*\s},2) -- (-0.5,1);
\draw[fill=gray!20] (0,0) -- ({0.5+\s},1) -- ({1.5-\s},1) -- (0,0);
\node[circle,inner sep=\is,outer sep=\os] (n1) at (-0.5,1) {\tiny$\bullet$};
\node[circle,inner sep=\is,outer sep=\os] (n0) at (0,0) {\tiny$\bullet$};
\node[circle,inner sep=0.7,outer sep=0.7] (n0b) at (0,0) {};
\draw[->] (n0b.-112) .. controls ({0-0.2},-0.3) and ({0+0.2},-0.3) .. (n0b.-68);
\draw 
(n1) edge (n0)
;
\node  at (-0.5,2+2*\s) {\scriptsize$8$};
\node  at (1,1+2*\s) {\scriptsize$6$};
\end{tikzpicture}
\end{array}
&\begin{array}{c}\times\end{array}&
%
%
\begin{array}{c}
\begin{tikzpicture}
\def\s{0.10}
\def\v{0.2}
\def\is{-1}
\def\os{-1}
\draw[fill=gray!20] (-0.5,1) -- ({-1+2*\s},2) -- ({0-2*\s},2) -- (-0.5,1);
\draw[fill=gray!20] (0,0) -- ({0.5+\s},1) -- ({1.5-\s},1) -- (0,0);
\node[circle,inner sep=\is,outer sep=\os] (n1) at (-0.5,1) {\tiny$\bullet$};
\node[circle,inner sep=\is,outer sep=\os] (n0) at (0,0) {\tiny$\bullet$};
\node[circle,inner sep=0.7,outer sep=0.7] (n0b) at (0,0) {};
\draw[->] (n0b.-112) .. controls ({0-0.2},-0.3) and ({0+0.2},-0.3) .. (n0b.-68);
\draw 
(n1) edge (n0)
;
\node  at (-0.5,2+2*\s) {\scriptsize$8$};
\node  at (1,1+2*\s) {\scriptsize$6$};
\end{tikzpicture}
\end{array}
\end{array}
\]
That $L_{1,2}\in\G(f)$ means that $f$ is isomorphic to the product of a $2$-nilpotent function in $F(1,16)$ and a constant function in $F(1,16)$, and indeed:
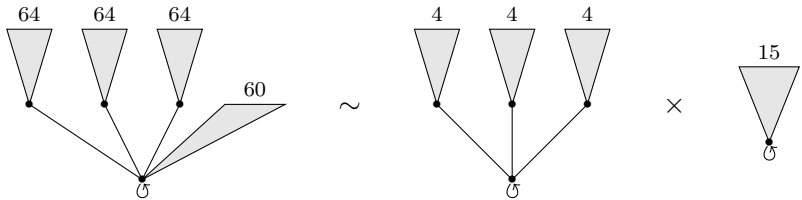
\[
\begin{array}{ccccc}
\begin{array}{c}
\begin{tikzpicture}
\def\s{0.10}
\def\v{0.2}
\def\is{-1}
\def\os{-1}
\draw[fill=gray!20] (-1.5,1) -- ({-2+2*\s},2) -- ({-1-2*\s},2) -- (-1.5,1);
\draw[fill=gray!20] (-0.5,1) -- ({-1+2*\s},2) -- ({0-2*\s},2) -- (-0.5,1);
\draw[fill=gray!20] (0.5,1) -- ({0+2*\s},2) -- ({1-2*\s},2) -- (0.5,1);
\draw[fill=gray!20] (0,0) -- ({1+\s},1) -- ({2-\s},1) -- (0,0);
\node[circle,inner sep=\is,outer sep=\os] (n1) at (-1.5,1) {\tiny$\bullet$};
\node[circle,inner sep=\is,outer sep=\os] (n2) at (-0.5,1) {\tiny$\bullet$};
\node[circle,inner sep=\is,outer sep=\os] (n3) at (0.5,1) {\tiny$\bullet$};
\node[circle,inner sep=\is,outer sep=\os] (n0) at (0,0) {\tiny$\bullet$};
\node[circle,inner sep=0.7,outer sep=0.7] (n0b) at (0,0) {};
\draw[->] (n0b.-112) .. controls ({0-0.2},-0.3) and ({0+0.2},-0.3) .. (n0b.-68);
\draw 
(n1) edge (n0)
(n2) edge (n0)
(n3) edge (n0)
;
\node  at (-1.5,2+2*\s) {\scriptsize$64$};
\node  at (-0.5,2+2*\s) {\scriptsize$64$};
\node  at (+0.5,2+2*\s) {\scriptsize$64$};
\node  at (1.5,1+2*\s) {\scriptsize$60$};
\end{tikzpicture}
\end{array}
&\begin{array}{c}\sim\end{array}&
%
\begin{array}{c}
\begin{tikzpicture}
\def\s{0.10}
\def\v{0.2}
\def\is{-1}
\def\os{-1}
\draw[fill=gray!20] (-1.5,1) -- ({-2+2*\s},2) -- ({-1-2*\s},2) -- (-1.5,1);
\draw[fill=gray!20] (-0.5,1) -- ({-1+2*\s},2) -- ({0-2*\s},2) -- (-0.5,1);
\draw[fill=gray!20] (0.5,1) -- ({0+2*\s},2) -- ({1-2*\s},2) -- (0.5,1);
\node[circle,inner sep=\is,outer sep=\os] (n1) at (-1.5,1) {\tiny$\bullet$};
\node[circle,inner sep=\is,outer sep=\os] (n2) at (-0.5,1) {\tiny$\bullet$};
\node[circle,inner sep=\is,outer sep=\os] (n3) at (0.5,1) {\tiny$\bullet$};
\node[circle,inner sep=\is,outer sep=\os] (n0) at (-0.5,0) {\tiny$\bullet$};
\node[circle,inner sep=0.7,outer sep=0.7] (n0b) at (-0.5,0) {};
\draw[->] (n0b.-112) .. controls ({-0.5-0.2},-0.3) and ({-0.5+0.2},-0.3) .. (n0b.-68);
\draw 
(n1) edge (n0)
(n2) edge (n0)
(n3) edge (n0)
;
\node  at (-1.5,2+2*\s) {\scriptsize$4$};
\node  at (-0.5,2+2*\s) {\scriptsize$4$};
\node  at (+0.5,2+2*\s) {\scriptsize$4$};
\end{tikzpicture}
\end{array}
&\begin{array}{c}\times\end{array}&
%
%
\begin{array}{c}
\begin{tikzpicture}
\def\s{0.10}
\def\v{0.2}
\def\is{-1}
\def\os{-1}
\draw[fill=gray!20] (0,0) -- ({-0.5+\s},1) -- ({0.5-\s},1) -- (0,0);
\node[circle,inner sep=\is,outer sep=\os] (n0) at (0,0) {\tiny$\bullet$};
\node[circle,inner sep=0.7,outer sep=0.7] (n0b) at (0,0) {};
\draw[->] (n0b.-112) .. controls ({0-0.2},-0.3) and ({0+0.2},-0.3) .. (n0b.-68);
\node  at (0,1+2*\s) {\scriptsize$15$};
\end{tikzpicture}
\end{array}
\end{array}
\]
\end{example}

\subsection{A characterization of universal functions}

In this section, we prove that the regular condition can be removed from Theorem \ref{thm:universal_regular_1}. We thus obtain the following characterization of universal functions in the general case (since any function which can be produced by an acyclic digraph of height $1$ is $2$-nilpotent (Theorem \ref{thm:robert}), this caracterization implies Theorem~\ref{thm:universal_1_intro} from the introduction). 

\begin{theorem}\label{thm:universal_1}
Let $n,q\geq 2$ and $f\in F(n,q)$. The following statements are equivalent:
\begin{itemize} 
\item
$f$ is universal,
\item
$f$ is $2$-nilpotent and $L_{1,n},L_{n,n}\in\G(f)$.
\end{itemize}
\end{theorem}

Since universal functions are $2$-nilpotent (Lemma~\ref{lem:2-nilpotent}), the first item trivially implies the second. For the converse direction, we assume that $f$ is $2$-nilpotent, and we proceed in three steps, using the following definition: for $m\in [n]$, let \EM{$L_{m,n}$} be the digraph on $[n]$ with a loop on each vertex $i\in [m]$ and no other arcs.  
\begin{itemize}
\item
First we give a very simple (and easy to check) necessary and sufficient condition for $\G(f)$ to contain $L_{1,n}$, which implies that the fixed point of $f$ has many pre-images (Lemma~\ref{lem:L1n}). 
\item
Second, we deduce from this characterization that $L_{1,n},L_{n,n}\in\G(f)$ implies $L_{m,n}\in\G(f)$ for all $m\in [n]$ (Lemma~\ref{lem:Lmn}). 
\item
Third, we prove that if $L_{m,n}\in\G(f)$ (which holds by Lemma~\ref{lem:Lmn}) and if the fixed point of $f$ has many pre-images (which holds by Lemma~\ref{lem:L1n}) then $\G(f)$ contains any digraph on $[n]$ with exactly $n-m$ sinks (Lemma~\ref{lem:n-m_sinks}). 
\end{itemize}
Together, this proves that $f$ is universal. 

\medskip
For the first step, we need a definition: we say that a $2$-nilpotent function $f\in F(n,q)$ has the \EM{pre-image property} if:
\begin{itemize}
\item $|f^{-1}(x)|$ is a multiple of $q^{n-1}$ for all $x\in\Im(f)$, and
\item the fixed point of $f$ has at least $|\Im(f)|q^{n-1}$ pre-images.
\end{itemize}
The characterization mentioned in the first step is the following. 

\begin{lemma}\label{lem:L1n}
Let $n,q\geq 2$ and let $f\in F(n,q)$ be a $2$-nilpotent function. The following statements are equivalent:
\begin{itemize}
\item $L_{1,n}\in\G(f)$,
\item $f$ has the pre-image property.
\end{itemize}
\end{lemma}

\begin{proof} Suppose, without loss, that $0^n$ is the fixed point of $f$, and let $r=|\Im(f)|$.

\medskip
Suppose first that $L_{1,n}\in \G(f)$ and let us prove that $f$ has the pre-image property. Since this property is invariant by isomorphism, we can suppose that $G(f)=L_{1,n}$. Since the local function $f_1$ only depends on input $1$, we regard $f_1$ as a transformation of $\Q{q}$. Then, since $f_i=\cst=0$ for all $1<i\leq n$, we have $f(x)=(f_1(x_1),0,\dots,0)$ for all $x\in\Q{q}^n$. Thus $\Im(f)=\Im(f_1)\times\{0\}^{n-1}$ and $f^{-1}(x)=f^{-1}_1(x_1)\times \Q{q}^{n-1}$ for all $x\in\Im(f)$. Thus $|f^{-1}(x)|$ is a multiple of $q^{n-1}$. Since $f$ is $2$-nilpotent, we deduce that $f_1$ is $2$-nilpotent and that its fixed point is $0$. Hence all the images of $f_1$ are contained in $f^{-1}_1(0)$, and thus $|f^{-1}_1(0)|\geq |\Im(f_1)|=r$ so that $|f^{-1}(0^n)|\geq  rq^{n-1}$. Thus $f$ has the pre-image property. 

\medskip
Conversely, suppose that $f$ has the pre-image property and let us prove that $L_{1,n}\in\G(f)$. Since $f$ has the pre-image property, for all $x\in\Im(f)$, there is a positive integer $\gamma(x)$ such that $|f^{-1}(x)|=\gamma(x)q^{n-1}$. Since the $r$ integers $\gamma(x)$ sum to $q$, we can assign to each $x\in\Im(f)$ a subset $A(x)$ of $\Q{q}$ of size $\gamma(x)$ such that $A(x)\cap A(y)=\emptyset$ for distinct $x,y\in\Im(f)$. There is then a permutation $\pi$ of $\Q{q}^n$ such that, for all $x\in\Im(f)$, 
\[
\pi(f^{-1}(x))=A(x)\times\Q{q}^{n-1}.
\]
Furthermore, since $f$ has the pre-image property, we have $rq^{n-1}\leq |f^{-1}(0^n)|=|A(0^n)|q^{n-1}$ and thus $r\leq |A(0^n)|$. Since $f$ is $2$-nilpotent, we have $\Im(f)\subseteq f^{-1}(0^n)$ and we deduce that the permutation $\pi$ can be chosen so that 
\[
\pi(\Im(f))\subseteq A(0^n)\times\{0\}^{n-1}.
\]
Let $h=\pi\circ f\circ \pi^{-1}$. Then $\Im(h)\subseteq A(0^n)\times\{0\}^{n-1}$ so $h_i=\cst=0$ for all $1<i\leq n$. Thus each $1<i\leq n$ is a source of $G(h)$. Furthermore, for any $x,y\in\Q{q}^n$ with $x_1=y_1$, we have  $x_1,y_1\in A(z)$ for some $z\in\Im(f)$. Thus $x,y\in A(z)\times\Q{q}^{n-1}=\pi(f^{-1}(z))$, that is, $f(\pi^{-1}(x))=f(\pi^{-1}(y))=z$ so that $h(x)=h(y)=\pi(z)$. So for any $x,y\in\Q{q}^n$ that only differ in some component $1< i\leq n$, we have $x_1=y_1$ and thus $h(x)=h(y)$. So each $1<i\leq n$ is an isolated vertex of $G(h)$, and since $h\neq\cst$ we deduce that $G(h)=L_{1,n}$. Thus $L_{1,n}\in \G(f)$. 
\end{proof}

Using the previous lemma, we obtain the following, which is the second step. 

\begin{lemma}\label{lem:Lmn}
Let $n,q\geq 2$ and let $f\in F(n,q)$ be a $2$-nilpotent function. The following statements are equivalent:
\begin{itemize} 
\item
$\G(f)$ contains $L_{1,n}$ and $L_{n,n}$,
\item
$\G(f)$ contains $L_{m,n}$ for all $m\in [n]$. 
\end{itemize}
\end{lemma}

\begin{proof}
It is clear that the second item implies the first. We prove the converse by induction on $n$. The case $n=2$ is trivial. So suppose that $n>2$. Let $f\in F(n,q)$ such that $\G(f)$ contains $L_{1,n}$ and $L_{n,n}$. Let $1<m<n$ and let us prove that $L_{m,n}\in\G(f)$.
 Suppose, without loss, that $0^n$ is the fixed point of~$f$. Given $x\in\Q{q}^{n-1}$ and $a\in\Q{q}$, we set $(x,a)=(x_1,\dots,x_{n-1},a)\in\Q{q}^n$.

\medskip
Since $G(f)=L_{n,n}$, each local function $f_i$ does not depend on input $j\neq i$. Hence we regard $f_i$ as a transformation of $\Q{q}$, so that $f(x)=(f_1(x_1),\dots,f_n(x_n))$ for all $x\in\Q{q}^n$. Let $f'\in F(n-1,q)$ defined by 
\[
f'(x)=(f_1(x_1),\dots,f_{n-1}(x_{n-1}))
\]
for all $x\in\Q{q}^{n-1}$. Since $f$ is $2$-nilpotent, $f'$ is $2$-nilpotent and its fixed point is $0^{n-1}$. Furthermore, $G(f')=L_{n-1,n-1}$ and, for all $x\in\Q{q}^{n-1}$ and $a\in\Q{q}$, we have 
\[
f'^{-1}(x)=f^{-1}_1(x_1)\times\dots\times f^{-1}_{n-1}(x_{n-1})
\]
and 
\[
f^{-1}(x,a)=f'^{-1}(x)\times f^{-1}_n(a).
\] 

\medskip
Let us prove that $L_{1,n-1}\in\G(f')$. By Lemma~\ref{lem:L1n}, $f$ has the pre-image property, that is, for all $x\in\Im(f)$ there is a positive integer $\gamma(x)$ such that $|f^{-1}(x)|=\gamma(x)q^{n-1}$ and $\gamma(0^n)\geq |\Im(f)|$. Given $x\in\Im(f')$ and  $a\in\Im(f_n)$, we have $(x,a)\in\Im(f)$ and we set 
\[
\gamma'(x)=\sum_{a\in\Im(f_n)}\gamma(x,a). 
\]
Then 
\[
\gamma'(x)q^{n-1}=\sum_{a\in\Im(f_n)}|f^{-1}(x,a)|=\sum_{a\in\Im(f_n)}|f'^{-1}(x)|\cdot|f_n^{-1}(a)|=q\cdot |f'^{-1}(x)|
\]
and thus $|f'^{-1}(x)|=\gamma'(x)q^{n-2}$. Since $\gamma'(0^{n-1})\geq \gamma(0^n)\geq |\Im(f)|\geq |\Im(f')|$ we deduce that $f'$ has the pre-image property. Hence, by Lemma~\ref{lem:L1n} we have $L_{1,n-1}\in\G(f')$. 

\medskip
Since $L_{1,n-1}$ and $L_{n-1,n-1} \in G(f')$, by induction hypothesis, there exists $h'\in F(L_{m-1,n-1},q)$ isomorphic to $f'$. Let $h\in F(n,q)$ defined by $h(x,a)=(h'(x),f_n(a))$ for all $x\in\Q{q}^{n-1}$ and $a\in\Q{q}$, which is isomorphic to $f$. Then $G(h)$ is obtained from $G(h')$ by adding vertex $n$ and a loop on it (since $f_n\neq \cst$). Thus $G(h)$ has $m$ loops (one on each vertex in $[m-1]\cup\{n\}$) and no other arcs. Hence $G(h)$ is isomorphic to $L_{m,n}$. Thus $L_{m,n}\in\G(f)$ because $G(h)\in\G(f)$ and $\G(f)$ is closed under isomorphism. 
\end{proof}

The last step is the next lemma, which is independent of the two previous ones.

\begin{lemma}\label{lem:n-m_sinks}
Let $n,q\geq 2$ and let $f\in F(n,q)$ be a $2$-nilpotent function. Let $1\leq m\leq n$ and suppose that the fixed point of $f$ has at least $|\Im(f)|(q-1)^{m-1}q^{n-m}$ pre-images. The following statements are equivalent:
\begin{itemize} 
\item
$\G(f)$ contains $L_{m,n}$,
\item
$\G(f)$ contains all the digraphs on $[n]$ with exactly $n-m$ sinks.
\end{itemize}
\end{lemma}

\begin{proof}
Since $L_{m,n}$ has exactly $n-m$ sinks, the second item implies the first. For the converse, suppose that $\G(f)$ contains $L_{m,n}$. We can suppose, without loss, that $G(f)=L_{m,n}$ and that the fixed point of $f$ is $0^n$. Then $f_i=\cst=0$ for all $m<i\leq n$. For all $i\in [m]$, since $f_i$ only depends on input $i$, we regard $f_i$ has a transformation of $\Q{q}$, so that $f(x)=(f_1(x_1),\dots,f_s(x_m),0,\dots,0)$ for all $x\in\Q{q}^n$. Consequently, 
\begin{equation}\label{eq:im}
0^n\in \Im(f)=\Im(f_1)\times\dots\times \Im(f_m)\times\{0\}^{n-m}.
\end{equation}
Furthermore, we have $f^{-1}(0^n)=f^{-1}_1(0)\times\dots\times f^{-1}_m(0)\times\Q{q}^{n-m}$ which is, by hypothesis, of size at least $r(q-1)^{m-1}q^{n-m}$, where $r=|\Im(f)|$. Since for all $i\in [m]$ we have $|f^{-1}_i(0)|\leq q-1$ (because $f_i$ is not constant), we deduce that $|f^{-1}_i(0)|(q-1)^{m-1}q^{n-m}\geq |f^{-1}(0^n)|\geq r(q-1)^{m-1}q^{n-m}$ and thus $|f^{-1}_i(0)|\geq r$. Hence we can assume, without loss, that $\Q{r}\subseteq f^{-1}_i(0)$ for all $i\in [n]$, so that 
\[
\Q{r}^n\subseteq f^{-1}(0^n). 
\]

\medskip
Let $G$ be any digraph on $[n]$ with exactly $n-m$ sinks. We will prove that $G\in\G(f)$. Since $\G(f)$ is closed under isomorphism, we can suppose that the set of sinks of $G$ is $]m,n]$. 

\medskip
For all $i\in [n]$, let $a_i\in\B^n$ defined by $a_{ij}=1$ if $G$ has an arc from $j$ to $i$, and $a_{ij}=0$ otherwise. For all $x\in\Q{q}^n$ we set $a_i\cdot f(x)=(a_{i1}\cdot f_1(x),\dots,a_{in}\cdot f_n(x))$; note that $a_i\cdot f(x)\in\Im(f)$ by \eqref{eq:im}. Let $\sigma$ be a bijection from $\Im(f)$ to $\Q{r}$ with $\sigma(0^n)=0$. Finally, let $h\in F(n,q)$ defined as follows: for all $x\in\Q{q}^n$ and $i\in [n]$, 
\[
h_i(x)=\sigma(a_i\cdot f(x)).
\] 
We will prove that $G(h)=G$ and then that $h\sim f$. 

\begin{itemize}
\item[(1)] $G(h)=G$.

Suppose that $G$ has an arc from $j$ to $i$. Then $j\in [m]$ so $f_j$ depends on input $j$: there are $x,y\in\Q{q}^n$ which only differ in $x_j\neq y_j$ such that $f_j(x)\neq f_j(y)$. Since $a_{ij}=1$, we have $a_i\cdot f(x)\neq a_i\cdot f(y)$, and since $\sigma$ is a permutation, we deduce that $h_i(x)\neq h_i(y)$. So $G(h)$ has an arc from $j$ to $i$. Conversely, suppose that $G(h)$ has an arc from $j$ to $i$. Then there are $x,y\in\Q{q}^n$ that only differ in $x_j\neq y_j$ such that $h_i(x)\neq h_i(y)$. Then $a_i\cdot f(x)\neq a_i\cdot f(y)$ so $f(x)\neq f(y)$ and since $x$ and $y$ only differ in $x_j\neq x_j$ and $G(f)=L_{m,n}$, we deduce that $f(x)$ and $f(y)$ only differ in $f_j(x)\neq f_j(y)$. Since $a_i\cdot f(x)\neq a_i\cdot f(y)$ this implies $a_{ij}=1$, and thus $G$ has an arc from $j$ to $i$. This proves (1).
\end{itemize}

\medskip
It remains to prove that $h\sim f$. 

\begin{itemize}
\item[(2)] {\em $h^2=0^n$.}

Let $x\in\Q{q}^n$. We have $h(x)\in \Im(\sigma)^n= \Q{r}^n\subseteq f^{-1}(0^n)$ and thus $f(h(x))=0^n$. We deduce that, for all $i\in [n]$ we have $h_i(h(x))=\sigma(a_i\cdot f(h(x)))=\sigma(a_i\cdot 0^n)=\sigma(0^n)=0$, and thus $h(h(x))=0^n$. This proves (2).
\end{itemize}

\begin{itemize}
\item[(3)] {\em For all $x,y\in \Q{q}^n$ we have $f(x)=f(y)$ if and only if $h(x)=h(y)$.}

If $f(x)=f(y)$ then it is clear that $h(x)=h(y)$. For the converse, suppose that $f(x)\neq f(y)$. Since $f_j=\cst=0$ for all $m<j\leq n$, we deduce that there is $j\in [m]$ with $f_j(x)\neq f_j(y)$. Since $f_j$ only depends on input $j$, we deduce that $x_j\neq y_j$. Since $j$ is not a sink of $G$, it has at least one out-neighbor in $G$, say $i$. So $a_{ij}=1$ and thus $a_i\cdot f(x)\neq a_i\cdot f(y)$. Since $\sigma$ is a permutation, we deduce that $h_i(x)\neq h_i(y)$. This proves (3).
\end{itemize}

Let $\pi:\Im(f)\to\Im(h)$ defined as follows: for all $x\in\Q{q}^n$, 
\[
\pi(f(x))=h(x);
\]
it is well defined since, by (3), if $f(x)=f(y)$ then $h(x)=h(y)$. It is clearly a surjection and, by (3), if $f(x)\neq f(y)$ then $\pi(f(x))=h(x)\neq h(y)= \pi(f(y))$. Hence $\pi$ is a bijection. Consequently, for all $x,y\in\Q{q}^n$ we have $f(y)=x$ if and only if $h(y)=\pi(f(y))=\pi(x)$. Hence, for all $x\in\Im(f)$ we have
\[
f^{-1}(x)=h^{-1}(\pi(x)).
\]
By (2), $h$ is $2$-nilpotent with fixed point $0^n$, and thus $\pi$ maps the fixed point of $f$ on the fixed point of $h$ since 
$\pi(0^n)=\pi(f(0^n))=h(0^n)=0^n$. Consequently, by Lemma~\ref{lem:2-nil_iso}, $h\sim f$. 
\end{proof}

We are now ready to prove Theorem~\ref{thm:universal_1}.

\begin{proof}[{\bf Proof of Theorem~\ref{thm:universal_1}}]
Since any universal function is $2$-nilpotent (Lemma \ref{lem:2-nilpotent}), the first item trivially implies the second. To prove the converse, suppose that $f$ is $2$-nilpotent and that $L_{1,n},L_{n,n}\in\G(f)$. Let $G$ be any digraph on $[n]$ with at least one arc. Then there exists an integer $m\in [n]$ such that $G$ has exactly $n-m$ sinks. By Lemma~\ref{lem:Lmn}, we have $L_{m,n}\in\G(f)$ and, by Lemma~\ref{lem:L1n}, the number of pre-images of the fixed point of $f$ is at least $|\Im(f)|q^{n-1}\geq |\Im(f)|(q-1)^{m-1}q^{n-m}$. Thus, by Lemma~\ref{lem:n-m_sinks}, we have $G\in\G(f)$. Hence $f$ is universal.  
\end{proof}

\begin{example}
Let $f\in F(2,9)$ be a function isomorphic to the following digraph (edges are oriented toward the loop):
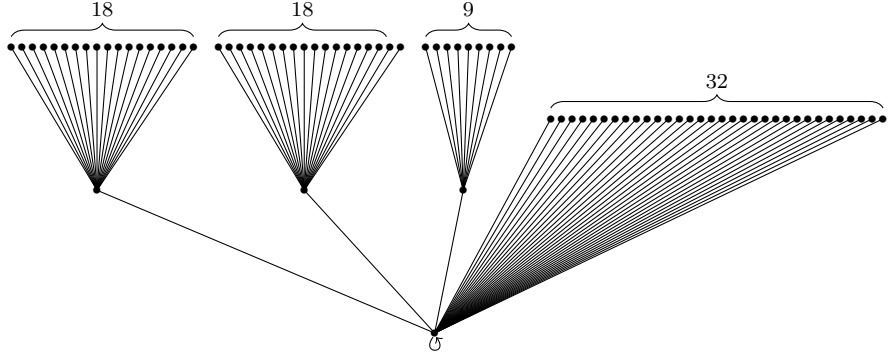
\[
\begin{tikzpicture}[scale=0.95]
\def\s{0.15}
\def\v{0.2}
\def\is{-1}
\def\os{-1}
\node[circle,inner sep=\is,outer sep=\os] (n1) at ({18*\s/2},2) {\tiny$\bullet$};
\node[circle,inner sep=\is,outer sep=\os] (n2) at ({\v+(18*\s)+18*\s/2},2) {\tiny$\bullet$};
\node[circle,inner sep=\is,outer sep=\os] (n3) at ({2*\v+(36*\s)+9*\s/2},2) {\tiny$\bullet$};
\node[circle,inner sep=\is,outer sep=\os] (n0) at ({(36*\s)+9*\s/2},0) {\tiny$\bullet$};
\node[circle,inner sep=0.7,outer sep=0.7] (n0b) at ({(36*\s)+9*\s/2},0) {};
\draw[->] (n0b.-112) .. controls ({(36*\s)+9*\s/2-0.2},-0.3) and ({(36*\s)+9*\s/2+0.2},-0.3) .. (n0b.-68);
\draw 
(n1) edge (n0)
(n2) edge (n0)
(n3) edge (n0)
;
\foreach \x in {1,...,18}{
	\node[circle,inner sep=\is,outer sep=-1] (\x) at ({0*\v+\x*\s},4) {\tiny$\bullet$};
	\draw (\x) edge (n1);
}
\foreach \x in {19,...,36}{
	\node[circle,inner sep=\is,outer sep=-1] (\x) at ({1*\v+\x*\s},4) {\tiny$\bullet$};
	\draw (\x) edge (n2);
}
\foreach \x in {37,...,45}{
	\node[circle,inner sep=\is,outer sep=-1] (\x) at ({2*\v+\x*\s},4) {\tiny$\bullet$};
	\draw (\x) edge (n3);
}
\foreach \x in {46,...,77}{
	\node[circle,inner sep=\is,outer sep=\os] (\x) at ({4*\v+\x*\s},3) {\tiny$\bullet$};
	\draw (\x) edge (n0);
}
\draw [decorate,decoration={brace,amplitude=5},xshift=0pt,yshift=0pt] 
({0*\v+1*\s},4+\s) -- ({0*\v+18*\s},4+\s) node [midway,above,outer sep=4] {\scriptsize{$18$}};
\draw [decorate,decoration={brace,amplitude=5},xshift=0pt,yshift=0pt] 
({1*\v+19*\s},4+\s) -- ({0*\v+36*\s},4+\s) node [midway,above,outer sep=4] {\scriptsize{$18$}};
\draw [decorate,decoration={brace,amplitude=5},xshift=0pt,yshift=0pt] 
({2*\v+37*\s},4+\s) -- ({2*\v+45*\s},4+\s) node [midway,above,outer sep=4] {\scriptsize{$9$}};
\draw [decorate,decoration={brace,amplitude=5},xshift=0pt,yshift=0pt] 
({4*\v+46*\s},3+\s) -- ({4*\v+77*\s},3+\s) node [midway,above,outer sep=4] {\scriptsize{$32$}};
\end{tikzpicture}
\]
Then $f$ is universal. This can be easily checked with Theorem~\ref{thm:universal_1}. For that, we have to show $f$ is $2$-nilpotent, which is obvious, and that $L_{1,2}$ and $L_{2,2}$ are in $\G(f)$. That $L_{2,2}\in\G(f)$ means that $f$  is isomorphic to the product of two $2$-nilpotent functions in $F(1,9)$, and indeed:
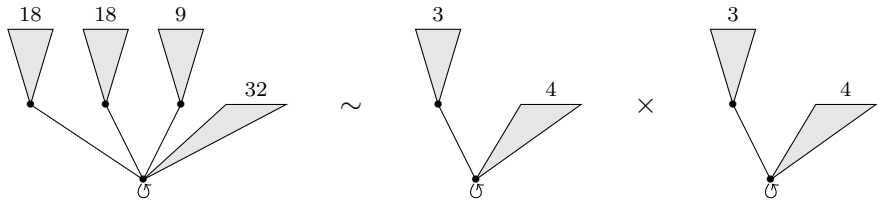
\[
\begin{array}{ccccc}
\begin{array}{c}
\begin{tikzpicture}
\def\s{0.10}
\def\v{0.2}
\def\is{-1}
\def\os{-1}
\draw[fill=gray!20] (-1.5,1) -- ({-2+2*\s},2) -- ({-1-2*\s},2) -- (-1.5,1);
\draw[fill=gray!20] (-0.5,1) -- ({-1+2*\s},2) -- ({0-2*\s},2) -- (-0.5,1);
\draw[fill=gray!20] (0.5,1) -- ({0+2*\s},2) -- ({1-2*\s},2) -- (0.5,1);
\draw[fill=gray!20] (0,0) -- ({1+\s},1) -- ({2-\s},1) -- (0,0);
\node[circle,inner sep=\is,outer sep=\os] (n1) at (-1.5,1) {\tiny$\bullet$};
\node[circle,inner sep=\is,outer sep=\os] (n2) at (-0.5,1) {\tiny$\bullet$};
\node[circle,inner sep=\is,outer sep=\os] (n3) at (0.5,1) {\tiny$\bullet$};
\node[circle,inner sep=\is,outer sep=\os] (n0) at (0,0) {\tiny$\bullet$};
\node[circle,inner sep=0.7,outer sep=0.7] (n0b) at (0,0) {};
\draw[->] (n0b.-112) .. controls ({0-0.2},-0.3) and ({0+0.2},-0.3) .. (n0b.-68);
\draw 
(n1) edge (n0)
(n2) edge (n0)
(n3) edge (n0)
;
\node  at (-1.5,2+2*\s) {\scriptsize$18$};
\node  at (-0.5,2+2*\s) {\scriptsize$18$};
\node  at (+0.5,2+2*\s) {\scriptsize$9$};
\node  at (1.5,1+2*\s) {\scriptsize$32$};
\end{tikzpicture}
\end{array}
&\begin{array}{c}\sim\end{array}&
%
\begin{array}{c}
\begin{tikzpicture}
\def\s{0.10}
\def\v{0.2}
\def\is{-1}
\def\os{-1}
\draw[fill=gray!20] (-0.5,1) -- ({-1+2*\s},2) -- ({0-2*\s},2) -- (-0.5,1);
\draw[fill=gray!20] (0,0) -- ({0.5+\s},1) -- ({1.5-\s},1) -- (0,0);
\node[circle,inner sep=\is,outer sep=\os] (n1) at (-0.5,1) {\tiny$\bullet$};
\node[circle,inner sep=\is,outer sep=\os] (n0) at (0,0) {\tiny$\bullet$};
\node[circle,inner sep=0.7,outer sep=0.7] (n0b) at (0,0) {};
\draw[->] (n0b.-112) .. controls ({0-0.2},-0.3) and ({0+0.2},-0.3) .. (n0b.-68);
\draw 
(n1) edge (n0)
;
\node  at (-0.5,2+2*\s) {\scriptsize$3$};
\node  at (1,1+2*\s) {\scriptsize$4$};
\end{tikzpicture}
\end{array}
&\begin{array}{c}\times\end{array}&
%
%
\begin{array}{c}
\begin{tikzpicture}
\def\s{0.10}
\def\v{0.2}
\def\is{-1}
\def\os{-1}
\draw[fill=gray!20] (-0.5,1) -- ({-1+2*\s},2) -- ({0-2*\s},2) -- (-0.5,1);
\draw[fill=gray!20] (0,0) -- ({0.5+\s},1) -- ({1.5-\s},1) -- (0,0);
\node[circle,inner sep=\is,outer sep=\os] (n1) at (-0.5,1) {\tiny$\bullet$};
\node[circle,inner sep=\is,outer sep=\os] (n0) at (0,0) {\tiny$\bullet$};
\node[circle,inner sep=0.7,outer sep=0.7] (n0b) at (0,0) {};
\draw[->] (n0b.-112) .. controls ({0-0.2},-0.3) and ({0+0.2},-0.3) .. (n0b.-68);
\draw 
(n1) edge (n0)
;
\node  at (-0.5,2+2*\s) {\scriptsize$3$};
\node  at (1,1+2*\s) {\scriptsize$4$};
\end{tikzpicture}
\end{array}
\end{array}
\]
That $L_{1,2}\in\G(f)$ means that $f$ is isomorphic to the product of a $2$-nilpotent function in $F(1,9)$ and a constant function in $F(1,9)$, and indeed:
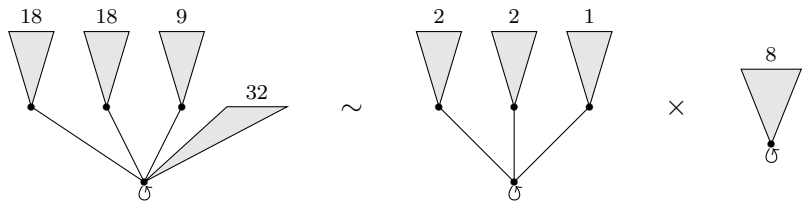
\[
\begin{array}{ccccc}
\begin{array}{c}
\begin{tikzpicture}
\def\s{0.10}
\def\v{0.2}
\def\is{-1}
\def\os{-1}
\draw[fill=gray!20] (-1.5,1) -- ({-2+2*\s},2) -- ({-1-2*\s},2) -- (-1.5,1);
\draw[fill=gray!20] (-0.5,1) -- ({-1+2*\s},2) -- ({0-2*\s},2) -- (-0.5,1);
\draw[fill=gray!20] (0.5,1) -- ({0+2*\s},2) -- ({1-2*\s},2) -- (0.5,1);
\draw[fill=gray!20] (0,0) -- ({1+\s},1) -- ({2-\s},1) -- (0,0);
\node[circle,inner sep=\is,outer sep=\os] (n1) at (-1.5,1) {\tiny$\bullet$};
\node[circle,inner sep=\is,outer sep=\os] (n2) at (-0.5,1) {\tiny$\bullet$};
\node[circle,inner sep=\is,outer sep=\os] (n3) at (0.5,1) {\tiny$\bullet$};
\node[circle,inner sep=\is,outer sep=\os] (n0) at (0,0) {\tiny$\bullet$};
\node[circle,inner sep=0.7,outer sep=0.7] (n0b) at (0,0) {};
\draw[->] (n0b.-112) .. controls ({0-0.2},-0.3) and ({0+0.2},-0.3) .. (n0b.-68);
\draw 
(n1) edge (n0)
(n2) edge (n0)
(n3) edge (n0)
;
\node  at (-1.5,2+2*\s) {\scriptsize$18$};
\node  at (-0.5,2+2*\s) {\scriptsize$18$};
\node  at (+0.5,2+2*\s) {\scriptsize$9$};
\node  at (1.5,1+2*\s) {\scriptsize$32$};
\end{tikzpicture}
\end{array}
&\begin{array}{c}\sim\end{array}&
%
\begin{array}{c}
\begin{tikzpicture}
\def\s{0.10}
\def\v{0.2}
\def\is{-1}
\def\os{-1}
\draw[fill=gray!20] (-1.5,1) -- ({-2+2*\s},2) -- ({-1-2*\s},2) -- (-1.5,1);
\draw[fill=gray!20] (-0.5,1) -- ({-1+2*\s},2) -- ({0-2*\s},2) -- (-0.5,1);
\draw[fill=gray!20] (0.5,1) -- ({0+2*\s},2) -- ({1-2*\s},2) -- (0.5,1);
\node[circle,inner sep=\is,outer sep=\os] (n1) at (-1.5,1) {\tiny$\bullet$};
\node[circle,inner sep=\is,outer sep=\os] (n2) at (-0.5,1) {\tiny$\bullet$};
\node[circle,inner sep=\is,outer sep=\os] (n3) at (0.5,1) {\tiny$\bullet$};
\node[circle,inner sep=\is,outer sep=\os] (n0) at (-0.5,0) {\tiny$\bullet$};
\node[circle,inner sep=0.7,outer sep=0.7] (n0b) at (-0.5,0) {};
\draw[->] (n0b.-112) .. controls ({-0.5-0.2},-0.3) and ({-0.5+0.2},-0.3) .. (n0b.-68);
\draw 
(n1) edge (n0)
(n2) edge (n0)
(n3) edge (n0)
;
\node  at (-1.5,2+2*\s) {\scriptsize$2$};
\node  at (-0.5,2+2*\s) {\scriptsize$2$};
\node  at (+0.5,2+2*\s) {\scriptsize$1$};
\end{tikzpicture}
\end{array}
&\begin{array}{c}\times\end{array}&
%
%
\begin{array}{c}
\begin{tikzpicture}
\def\s{0.10}
\def\v{0.2}
\def\is{-1}
\def\os{-1}
\draw[fill=gray!20] (0,0) -- ({-0.5+\s},1) -- ({0.5-\s},1) -- (0,0);
\node[circle,inner sep=\is,outer sep=\os] (n0) at (0,0) {\tiny$\bullet$};
\node[circle,inner sep=0.7,outer sep=0.7] (n0b) at (0,0) {};
\draw[->] (n0b.-112) .. controls ({0-0.2},-0.3) and ({0+0.2},-0.3) .. (n0b.-68);
\node  at (0,1+2*\s) {\scriptsize$8$};
\end{tikzpicture}
\end{array}
\end{array}
\]
\end{example}

\subsection{Characterization of sets \BF{$F(n,q)$} containing a universal function}

In this section, we give the following characterization of the couples $(n,q)$ such that $F(n,q)$ contains a universal function (this is a restatement of Theorem~\ref{thm:universal_2_intro} from the introduction). It shows that $q$ must have at least $n$ prime factors (with multiplicity). 

\begin{theorem}\label{thm:universal_2}
For $n,q\geq 2$, the following statements are equivalent:
\begin{itemize} 
\item
$F(n,q)$ contains a universal function,
\item
$q$ can be expressed as the product of $n$ integers $q=q_1q_2\dots q_n$ with $\prod_{i=1}^n(q_i-1)\geq 2^n$. 
\end{itemize}
\end{theorem}

We divide the proof in two steps. The first is a simplification step:  whenever $F(n,q)$ contains a universal function, it can be chosen to have rank $2^n$ (which is the minimal rank of any universal function with $n$ components).  

\begin{lemma}\label{lem:rank2n}
For $n,q\geq 2$, the following statements are equivalent:
\begin{itemize} 
\item
$F(n,q)$ contains a universal function,
\item
$F(n,q)$ contains a universal function of rank $2^n$.
\end{itemize}
\end{lemma}

\begin{proof}
The second item trivially implies the first. To prove the converse, suppose that $f\in F(n,q)$ is universal and let us prove that $F(n,q)$ contains a universal function of rank $2^n$. Suppose, without loss, that $G(f)=L_{n,n}$. Since $f_i$ only depends on $i$, we regard $f_i$ as a transformation of $\Q{q}$, so that $f(x)=(f_1(x_1),\dots,f_n(x_n))$ for all $x\in\Q{q}^n$. For all $i\in [n]$, we have $f_i\neq\cst$ so we can suppose, without loss, that $\B\subseteq\Im(f_i)$, so that $\B^n\subseteq \Im(f)$. By Lemma~\ref{lem:2-nilpotent}, $f$ is $2$-nilpotent, and we can suppose, without loss, that $0^n$ is the fixed point of $f$.  

\medskip
Let $\tilde f\in F(n,q)$ defined as follows: for all $x\in\Q{q}^n$ and $i\in [n]$, 
\[
\tilde f_i(x)=\tilde f_i(x_i)=\min(f_i(x_i),1). 
\]
We will prove that $\tilde f$ is universal of rank $2^n$. Since $f_i\neq\cst$ and $f_i(0)=0$, we have $\tilde f_i\neq\cst$, thus  $G(\tilde f)=L_{n,n}$ and $\Im(\tilde f)=\B^n\subseteq \Im(f)$. So $\tilde f$ is of rank $2^n$, and since $f$ is $2$-nilpotent, we deduce that $f(\tilde f(x))=0$, and thus $\tilde f^2(x)=0$, so $\tilde f$ is also $2$-nilpotent with fixed point $0^n$. By Theorem~\ref{thm:universal_1} and Lemma~\ref{lem:L1n}, it remains to prove that $\tilde f$ has the pre-image property. 

\medskip
Fix any $x\in\B^n$ and let $X=X_1\times\cdots\times X_n$ where $X_i=\{0\}$ if $x_i=0$ and $X_i=\Im(f_i)\setminus\{0\}$ otherwise. For any $y\in\Q{q}^n$ we have $\tilde f_i(y_i)=0$ iff $f_i(y_i)=0$ and thus $\tilde f_i(y_i)=x_i$ iff $f_i(y_i)\in X_i$. We deduce that $\tilde f(y)=x$ iff $f(y)\in X$ and thus $f^{-1}(x)=f^{-1}(X)$. So
\[
|\tilde f^{-1}(x)|=|f^{-1}(X)|=\sum_{z\in X} |f^{-1}(z)|.  
\]
By Lemma~\ref{lem:L1n}, $f$ has the pre-image property, thus each $|f^{-1}(z)|$ if a multiple of $q^{n-1}$ and we deduce that $|f^{-1}(x)|$ is a multiple of $q^{n-1}$. For $x=0^n$, we have $X=\{0^n\}$ and since $f$ has the pre-image property we obtain 
\[
|\tilde f^{-1}(0^n)|=|f^{-1}(0^n)|\geq |\Im(f)|q^{n-1}\geq|\Im(\tilde f)|q^{n-1}.  
\]
So $\tilde f$ has the pre-image property.
\end{proof}

Hence, to prove Theorem~\ref{thm:universal_2} it is thus sufficient to prove that $F(n,q)$ contains a universal function of rank $2^n$ if and only if the second item of Theorem~\ref{thm:universal_2}  is satisfied.

\begin{lemma}\label{lem:2n_q}
For $n,q\geq 2$, the following statements are equivalent:
\begin{itemize} 
\item
$F(n,q)$ contains a universal function of rank $2^n$,
\item
$q$ can be expressed as the product of $n$ integers $q=q_1q_2\dots q_n$ with $\prod_{i=1}^n(q_i-1)\geq 2^n$. 
\end{itemize}
\end{lemma}

\begin{proof}
We first prove that the second item implies the first. Actually, we prove the following strengthening:
\begin{itemize}
\item[(1)] {\em Let $q=q_1q_2\dots q_n$ be the product of $n$ integers and let $\alpha:=\prod_{i=1}^n(q_i-1)$. If $\alpha\geq 2^n$ then $F(n,q)$ contains a universal function of rank $2^n$ whose fixed point has exactly $\alpha q^{n-1}$ pre-images.}

For all $i\in [n]$, let 
\[
a^i_0=(q_i-1)\frac{q}{q_i},\qquad a^i_1=\frac{q}{q_i}.
\]
We have $a^i_0+a^i_1=q$, and also $a^i_0\geq 2$: this is obvious if $q_i\geq 3$ and otherwise, since $\alpha\geq 2^n$, we have $q_i=2$ and thus $a^i_0=q/q_i\geq \alpha/q_i\geq 2^{n-1}\geq 2$. Let $f_i:\Q{q}\to\{0,1\}$ defined by: for all $\ell\in\Q{q}$,  \[
f_i(\ell)=
\left\{\begin{array}{ll}
0&\textrm{if $\ell<a^i_0$,}\\
1&\textrm{otherwise.}
\end{array}
\right.
\]
Since $a^i_0\geq 2$ we have $f_i(1)=f_i(0)=0$ and thus $f_i$ is $2$-nilpotent. Let $f\in F(n,q)$ defined by $f(x)=(f_1(x_1),\dots,f_n(x_n))$ for all $x\in\Q{q}^n$. Since each $f_i$ is $2$-nilpotent, $f$ is 2-nilpotent and its fixed point is $0^n$. Furthermore, $\Im(f)=\B^n$ and $G(f)=L_{n,n}$. Also, for all $x\in\B^n$, we have  
\[
|f^{-1}(x)|=\prod_{i\in [n]}|f^{-1}_i(x_i)|=a^1_{x_1}a^2_{x_2}\cdots a^n_{x_n}.
\] 
Since $a^i_{x^i}$ is a multiple of $q/q_i$, we deduce that $|f^{-1}(x)|$ is a multiple of $\prod_{i=1}^n q/q_i=q^{n-1}$. 
Furthermore, $|f^{-1}(0^n)|=a^1_0a^2_0\cdots a^n_0= \alpha q^{n-1}\geq  2^nq^{n-1}=|\Im(f)|q^{n-1}$, and thus $f$ has the pre-image property. We deduce from Lemma~\ref{lem:L1n} that $L_{1,n}\in \G(f)$, and then from Theorem~\ref{thm:universal_1} that $f$ is universal. This proves (1). 
\end{itemize}

\medskip
We now prove that the first item implies the second. So suppose that $F(n,q)$ contains a universal function $f$ of rank $2^n$, and suppose that, for these properties, the number of pre-images of the fixed point is maximal. Suppose, without loss, that $G(f)=L_{n,n}$. Since $f_i$ only depends on $i$, we regard $f_i$ as a transformation of $\Q{q}$, so that $f(x)=(f_1(x_1),\dots,f_n(x_n))$ for all $x\in\Q{q}^n$. Since $f$ has rank $2^n$, we can suppose, without loss, that $\Im(f)=\B^n$. By Lemma~\ref{lem:2-nilpotent}, $f$ is $2$-nilpotent, and we can suppose, without loss, that $0^n$ is the fixed point of $f$. To simplify notations, we set 
\[
a^i_0=|f^{-1}_i(0)|,\qquad a^i_1=|f^{-1}_i(1)|, 
\]
so that $|f^{-1}(x)|=a^1_{x_1}a^2_{x_2}\cdots a^n_{x_n}$. Then, for all $i\in [n]$, we set 
\[
q_i=\frac{q}{\gcd(a^i_0,q)},
\]
which is an integer exceeding $1$ since $a^i_0<q$.  

\begin{itemize}
\item[(2)] {\em $q_1q_2\dots q_n$ divides $q$.}

For any prime $p$ and positive integer $m$, let $\nu_p(m)$ be the $p$-adic order of $m$, that is, the maximum integer $\alpha\geq 0$ such that $p^\alpha$ divides $m$. Then $q_1q_2\dots q_n$ divides $q$ if and only if, for any prime $p$, the sum of the $\nu_p(q_i)$ is at most $\nu_p(q)$. Suppose, for a contradiction, that 
\[
\sum_{i\in[n]} \nu_p(q_i)> \nu_p(q)
\]
for some prime $p$. Let $I$ be the set of $i\in [n]$ such that $\nu_p(a^i_0)<\nu_p(q)$. Since $\nu_p(q_i)=\nu_p(q)-\min(\nu_p(a^i_0),\nu_p(q))$, we have $\nu_p(q_i)=\nu_p(q)-\nu_p(a^i_0)$ if $i\in I$ and $\nu_p(q_i)=0$ otherwise. Thus 
\[
\sum_{i\in I}^n \nu_p(q)-\nu_p(a^i_0)=\sum_{i\in [n]} \nu_p(q_i)>\nu_p(q),
\]
so
\[
\sum_{i\in I}\nu_p(a^i_0)<(|I|-1)\nu_p(q).
\] 
Since $q=a^i_0+a^i_1$ we have $\nu_p(a^i_0) \leq \nu_p(q)$ or $\nu_p(a^i_1) \leq \nu_p(q)$. Thus there exists $\ell_i\in\B$ such that $\nu_p(a^i_{\ell_i})\leq\nu_p(q)$. Let $x\in \B^n$ defined by $x_i=0$ for $i\in I$ and $x_i=\ell_i$ for $i\in [n]\setminus I$. We have 
\[
\sum_{i\in [n]}\nu_p(a^i_{x_i})
=\sum_{i\in I}\nu_p(a^i_0)+\sum_{i\in [n]\setminus I}\nu_p(a^i_{\ell_i})
<(|I|-1)\nu_p(q)+(n-|I|)\nu_p(q)=(n-1)\nu_p(q).
\] 
But since $f$ has the pre-image property, $a^1_{x_1}a^2_{x_2}\dots a^{n}_{x_n}$ is a multiple of $q^{n-1}$, so 
\[
\sum_{i\in [n]}\nu_p(a^i_{x_i})\geq (n-1)\nu_p(q),
\]
a contradiction. This proves (2).
\end{itemize}

\begin{itemize}
\item[(3)] {\em $a^i_0\leq (q_i-1)q/q_i$ for all $i\in [n]$.}

Indeed, 
\[
\frac{\lcm(a^i_0,q)}{q}\cdot \frac{q}{q_i}
=\frac{\lcm(a^i_0,q)}{q}\cdot\frac{q\cdot \gcd(a^i_0,q)}{q}
=\frac{\lcm(a^i_0,q)\cdot \gcd(a^i_0,q)}{q}
=\frac{a^i_0\cdot q}{q}=a^i_0
\]
and since $\lcm(a^i_0,q)/q$ is an integer, and $\lcm(a^i_0,q)/q<q_i$ (since otherwise $a^i_0\geq q$), we have $\lcm(a^i_0,q)/q\leq q_i-1$ and thus $a^i_0\leq (q_i-1)q/q_i$. This proves (3).
\end{itemize}

We are now in position to conclude. Let $\beta=q/(q_1q_2\dots q_n)$, which is a positive integer by~(2). Since $f$ has the pre-image property and by (3) we have 
\begin{eqnarray*}
2^nq^{n-1}\leq |f^{-1}(0^n)|
&=&a^1_0a^2_0\dots a^n_0\\
&\leq& \frac{q^n}{q_1q_2\dots q_n}\prod_{i=1}^n(q_i-1)\\
&=& q^{n-1}\beta\prod_{i=1}^n(q_i-1)\\
&=&q^{n-1}\beta(q_1-1)\prod_{i=2}^n(q_i-1)\\
&\leq &q^{n-1}(\beta q_1-1)\prod_{i=2}^n(q_i-1).
\end{eqnarray*}
So $(\beta q_1-1)\prod_{i=2}^n(q_i-1)\geq 2^n$. Since $q=(\beta q_1)q_2\cdots q_n$, by (1) there exists a universal function $\tilde f\in F(n,q)$ of rank $2^n$ whose fixed point, say $0^n$ without loss, has exactly $q^{n-1}(\beta q_1-1)\prod_{i=2}^n(q_i-1)$ pre-images. But if $\beta\geq 2$ then $\beta(q_1-1)<(\beta q_1-1)$ thus $|f^{-1}(0^n)|<|\tilde f^{-1}(0^n)|$ and this contradicts the choice of $f$. Consequently, $\beta=1$, that is $q=q_1q_2\dots q_n$ and we deduce that $\prod_{i=1}^n(q_i-1)\geq 2^n$.  
\end{proof}

By, Theorem~\ref{thm:universal_2}, $F(n,q)$ contains a universal function for $q=3^n$. We prove below that, for $n$ fixed, $q=3^n$ is actually the smallest alphabet size which allows the existence of a universal function. 

\begin{corollary}
For $n,q\geq 2$, if $F(n,q)$ contains a universal function then $q\geq 3^n$. 
\end{corollary}

\begin{proof}
Suppose that $F(n,q)$ contains a universal function. By Theorem~\ref{thm:universal_2}, $q$ can be expressed as the product of $n$ integers $q=q_1q_2\dots q_n$ with $\prod_{i=1}^n(q_i-1)\geq 2^n$. For all $i\in [n]$, let $a_i=(q_i-1)q/q_i$ and $b_i=q/q_i$. By the inequality of arithmetic and geometric means, we have 
\[
\left( \frac{1}{n} \sum_{i = 1}^{n} a_i \right)^n \geq a_1 a_1 \cdots a_n =q^{n-1}\prod_{i=1}^n(q_i-1)\geq 2^nq^{n-1}
\]
and
\[
\left(\frac{1}{n} \sum_{i = 1}^{n} b_i \right)^n \geq b_1b_2 \cdots b_n = q^{n-1}.
\] 
Since $a_i+b_i= q$, we deduce that
\[ 
q = \frac{1}{n}\sum_{i = 1}^{n} (a_i + b_i)=\left(\frac{1}{n}\sum_{i = 1}^{n} a_i\right)+\left(\frac{1}{n}\sum_{i = 1}^{n} b_i\right) \geq 2q^{\frac{n-1}{n}}+q^{\frac{n-1}{n}}=3q^{\frac{n-1}{n}}.
\]
So $q^n \geq 3^nq^{n-1}$ and thus $q\geq 3^n$. 
\end{proof}

\section{Almost universal functions}\label{sec:almost_universal}

For $n\geq 1$ and $q\geq 2$, let 
\[
\gamma_q(n)=\max_{f\in F(n,q)}|\G(f)|.
\]
For instance, $\gamma_q(n)=2^{n^2}-1$ if and only if $F(n,q)$ contains a universal function. For that, we saw that $q$ must be exponential with $n$. But what happens for $q$ fixed when $n$ is large, which is the usual situation. Let us say there is a \EM{$q$-almost universal function} if 
\[
\gamma_q(n)/2^{n^2}\to 1\textrm{ as }n\to\infty,
\]
that is, if there is $f\in F(n,q)$ such that the probability that a random digraph on $[n]$ belongs to $\G(f)$ tends to $1$ as $n$ tends to $\infty$. In this section, we prove the following. 

\begin{theorem}\label{thm:almost_universal}
There exists a $q$-almost universal function for all $q\geq 3$.
\end{theorem}

First, we exhibit, for each $n\geq 1$ and $q\geq 3$, a function $f\in F(n,q)$, with $2^n$ periodic points, such that $\G(f)$ contains all the Hamiltonian digraphs on $[n]$. This is Theorem \ref{thm:Hamiltonian1} stated in the introduction. This proves Theorem~\ref{thm:almost_universal} since it is well known that almost all digraphs are Hamiltonian \cite{W73} (almost all digraphs have some property $P$ if, among the digraphs on $[n]$, the fraction of those that satisfy $P$ tends to $1$ as $n\to\infty$). 

\medskip
Second, for many alphabet sizes $q$, we exhibit a $2$-nilpotent function $f\in F(n,q)$ such that $\G(f)$ contains much more than the Hamiltonian digraphs on $[n]$. We then show that this construction gives an alternative proof that a $q$-almost universal function exists for all $q\geq 4$. 

\subsection{Hamiltonian digraphs}

In this subsection, we prove Theorem \ref{thm:Hamiltonian1}, that we restate. 

\begin{theorem}\label{thm:Hamiltonian2}
For all $n\geq 1$ and $q\geq 3$, there exists $f\in F(n,q)$ with $2^n$ periodic points, such that $\G(f)$ contains all the Hamiltonian digraphs on $[n]$.
\end{theorem}

\begin{proof}
Let $\ell\geq 2$ be an integer. We first construct a function $f\in F(n,2\ell)$ with the following two properties: $f$ has $2^n$ periodic points, and $\G(f)$ contains all the Hamiltonian digraphs on $[n]$. We then consider a restriction $f'$ of $f$, which belongs to $F(n,2\ell-1)$, and we prove that $f'$ has the same two properties. Together this proves the theorem. 

\medskip
The construction is based on the cyclic shift $\sigma\in F(n,2)$ defined by 
\[
\sigma(x_1,x_2,\dots,x_n)=(x_n,x_1,\dots,x_{n-1})
\]
for all $x\in\B^n$. This is a permutation whose interaction graph is $C_n$. Then, identifying $\Q{2\ell}^n$ with $\B^n\times\Q{\ell}^n$, we define $f\in F(n,2\ell)$ by 
\[
f(x,y)=(\sigma(x),0^n)
\]
for all $x\in\B^n$ and $y\in\Q{\ell}^n$. Since $\sigma$ is a permutation, it is clear that $f$ has $2^n$ images, which are all periodic. Moreover, $G(f)=G(\sigma)=C_n$. We will now prove that $\G(f)$ contains all the Hamiltonian digraphs on $[n]$. 

\medskip
So let $H$ be any hamiltonian digraph on $[n]$. Since $H$ is Hamiltonian, it is isomorphic to a digraph $G$ on $[n]$ containing $C_n$. Since $\G(f)$ is closed under isomorphism, we have $H\in\G(f)$ if and only if $G\in\G(f)$, and it is more convenient to prove that $G\in\G(f)$. Given $i\in [n]$, we denote by $\phi(i)$ the in-neighbor of $i$ in 
$C_n$, and we denote by $M_i$ the in-neighbors of $i$ in $G$ distinct from $\phi(i)$. 

\medskip
Let $g\in F(n,2\ell)$ defined by 
\[
g(x,y)=(\sigma(x),h(x))
\]
where $h$ is the function in $F(n,\ell)$ defined as follows: for all $i\in [n]$ and $x\in\Q{\ell}^n$, 
\[
h_i(x)=
\left\{
\begin{array}{ll}
1&\textrm{if $x_{\phi(i)}=0$ and $x_j=1$ for all $j\in M_i$},\\
0&\textrm{otherwise}.
\end{array}
\right.
\]

\medskip
We easily check that $G(h)=G$. Actually, since $\Im(h)=\B^n$, the restriction $h''$ of $h$ on $\B^n$ belongs to $F(n,2)$ and we easily check the stronger property that $G(h'')=G$. Consequently, 
\[
G(g)=G(\sigma)\cup G(h'')=C_n\cup G=G,
\]
thus it only remains to prove that $g\sim f$. 

\medskip
Let $\pi\in F(n,2\ell)$ defined by:
\[
\pi(x,y)=
\left\{
\begin{array}{ll}
(x,y\oplus h(\sigma^{-1}(x)))&\textrm{if $y\in\B^n$},\\
(x,y)&\textrm{otherwise}.
\end{array}
\right.
\]
We easily check that $\pi$ is a permutation of $\B^n\times\Q{\ell}^n$. Furthermore, for all $(x,y)\in\B^n\times\Q{\ell}^n$, noting that $h(x)\in\B^n$ (to obtain the second equality), and noting that $\pi(x,y)=(x,z)$ for some configuration $z\in\Q{\ell}^n$ (to obtain the last equality), we have  
\begin{eqnarray*}
\pi(g(x,y))
&=&\pi(\sigma(x),h(x))\\
&=&(\sigma(x),h(x)\oplus h(\sigma^{-1}(\sigma(x))))\\
&=&(\sigma(x),h(x)\oplus h(x))\\
&=&(\sigma(x),0^n)\\
&=&f(x,z)\\
&=&f(\pi(x,y)).
\end{eqnarray*}
So $g\sim f$ and thus $G\in\G(f)$. 

\medskip
We now construct, from $f$, a function $f'\in F(n,2\ell-1)$ with the same properties: $f'$ has $2^n$ periodic points and $\G(f')$ contains all the Hamiltonian digraphs on $[n]$. We identify the set of configurations $\Q{2\ell-1}^n$ with the set $\Lambda$ of configurations $(x,y)\in\B^n\times\Q{\ell}^n$ such that $(x_i,y_i)\neq (1,\ell-1)$ for all $i\in[n]$. 

\medskip
Let $f'$ be the restriction of $f$ on $\Lambda$. Since $\Im(f)=\B^n\times\{0^n\}\subseteq \Lambda$, we have $f'\in F(n,2\ell-1)$ and $f'$ and $f$ have the same periodic points. It remains to prove that $\G(f')$ contains all the Hamiltonian digraphs on $[n]$ and, as argued above, it is sufficient to prove that $G\in \G(f')$. 

\medskip
Let $g'$ be the restriction of $g$ on $\Lambda$. For all $(x,y)\in\Lambda$, if $h_i(x)<\ell-1$, then it is clear that $g(x,y)\in\Lambda$. Otherwise, $\ell=2$ and $h_i(x)=1$. Thus $x_{\phi(i)}=0$ and thus $\sigma_i(x)=x_{\phi(i)}=0$, so that $g_i(x,y)=(0,1) \neq (1,\ell-1)$. We deduce that $\Im(g)\subseteq \Lambda$, and thus $g'\in F(n,2\ell-1)$. We obviously have $G(g')\subseteq G(g)=G$ and to prove the converse inclusion, suppose that $G$ has an arc from $j$ to $i$. Since $G(h'')=G$, there exists $x\in\B^n$, which only differ in $x_j\neq y_j$ such that $h_i(x)\neq h_i(y)$. Since $x'=(x,0^n)$ and $y'=(y,0^n)$ belongs to $\Lambda$, we have $g'_i(x')=g_i(x')=(\sigma_i(x),h_i(x))\neq (\sigma_i(y),h_i(y))=g_i(y')=g'_i(y)$ and since $x'$ and $y'$ only differ in $x'_j\neq y'_j$ we deduce that $G(g')$ has an arc from $j$ to $i$. Thus $G(g')=G$. 

\medskip
Hence to prove that $G\in\G(f')$, it remains to prove that $g'\sim f'$. We first prove that $\pi(\Lambda)\subseteq\Lambda$. So let $(x,y)\in\Lambda$. If $y\not\in\B^n$ then $\pi(x,y)=(x,y)\in\Lambda$. So suppose that $y\in\B^n$ and let $z=\sigma^{-1}(x)$. We have to prove that $\pi_i(x,y)=(x_i,y_i\oplus h_i(z))\neq (1,\ell-1)$. This is obvious if $x_i=0$ or $l>2$, so suppose that $x_i=1$ and $\ell=2$, which implies $y_i=0$ since $(x,y)\in\Lambda$. Then $z_{\phi(i)}=x_i=1$ so $h_i(z)=0$ and we deduce that $\pi_i(x,y)=(1,0\oplus 0)=(1,0) \neq (1,l-1)$. So $\pi(\Lambda)\subseteq\Lambda$ and we deduce that the restriction $\pi'$ of $\pi$ on $\Lambda$ is a permutation of $\Lambda$. Since $\pi$ is an isomorphism between $f$ and $g$, we deduce that $\pi'$ is an isomorphism between $f'$ and $g'$. So $g'\sim f'$ and thus $G\in\G(f')$. 
\end{proof}

\subsection{Permutable digraphs}

Given $d\geq 2$, a digraph $G$ on $[n]$ is \EM{$d$-permutable} if $F(G,d)$ contains a permutation. In this subsection, we prove that for many alphabet sizes $q$, there exists a function in $F(n,q)$ which can be produced by all the $d$-permutable digraphs on $[n]$, for some divisors $d$ of $q$ or $q+1$.  

\begin{theorem}\label{thm:d-permutable}
For all $n\geq 1$, $q\geq 4$ and divisor $d$ of $q$ or $q+1$ satisfying $4\leq d^2\leq q$, there exists a $2$-nilpotent function $f\in F(n,q)$ such that $\G(f)$ contains all the $d$-permutable digraphs on $[n]$.   
\end{theorem}

This is interesting since there are many permutable digraphs. Indeed, let us say that $G$ is \EM{coverable} if the vertices of $G$ can be spanned with vertex-disjoint cycles. Gadouleau \cite{G18rank} shows that ``coverable'' and ``$q$-permutable'' are closely related properties.

\begin{theorem}[\cite{G18rank}]\label{thm:G18}
Let $n\geq 1$, $q\geq 2$ and let $G$ be a digraph on $[n]$. If $G$ is $q$-permutable, then $G$ is coverable, and the converse holds for $q\geq 3$ but fails for $q=2$.
\end{theorem}

For a given alphabet size $q\geq 4$, let $d$ be the largest divisor of $q$ or $q+1$ such that $4\leq d^2\leq q$; it exists since one of $q,q+1$ is even. If $d\geq 3$ then, by Theorems~\ref{thm:d-permutable} and \ref{thm:G18}, there exists $f\in F(n,q)$ such that $\G(f)$ contains all the coverable digraphs on $[n]$: this improves Theorem~\ref{thm:Hamiltonian2} since Hamiltonian digraphs are coverable. But even for very large $q$, we can have $d=2$. It is for instance the case if $q$ and $(q+1)/2$ are primes, and that an infinity of such $q$ exists is a difficult open problem (related to Dickson's conjecture). Actually, for $q\neq 7,9,17$, we have $d=2$ if and only if one of the following two conditions is true: (1) $p$ and $(p+1)/2$ are primes; (2) $p/2$ is prime and $p+1$ is prime or the square of a prime. 
%
%
However, if $d=2$ then, by Theorem~\ref{thm:d-permutable}, there exists $f\in F(n,q)$ such that $\G(f)$ contains all the $2$-permutable digraphs on $[n]$; this is still a large family of coverable digraphs since we prove the following (this is a restatement of Lemma~\ref{lem:2-permutable_intro} from the introduction):

\begin{lemma}\label{lem:2-permutable}
Almost all coverable digraphs are $2$-permutable. 
\end{lemma}

This completes Gadouleau's theorem, and since almost all digraphs are coverable, this gives an alternative proof of Theorem~\ref{thm:almost_universal} for $q\geq 4$. 

\medskip
We now prove Theorem~\ref{thm:d-permutable} and Lemma~\ref{lem:2-permutable}.

\begin{proof}[\BF{Proof of Theorem \ref{thm:d-permutable}}]
Since $d$ divides $q$ or $q+1$ and $4\leq d^2\leq q$, there exists a positive integer $\ell$ such that either $q=\ell d$ and $d\leq\ell$, or $q=\ell d-1$ and $d<\ell$. So it is sufficient to prove that, given $2\leq d\leq\ell$, $F(n,\ell d)$ contains a $2$-nilpotent function which can be generated by all the $d$-permutable digraphs on $[n]$, and if $d<\ell$ then $F(n,\ell d-1)$ also contains such a function. 

\medskip
We start with $F(n,\ell d)$. We identify the set of configurations $\Q{\ell d}^n$ with $\Q{\ell}^n\times\Q{d}^n$. Hence a configuration is a couple $(x,y)$ with $x\in\Q{\ell}^n$ and $y\in\Q{d}^n$. 

\medskip
Let $\sigma$ be the function from $\Q{\ell}^n$ to $\Q{ d}^n$ defined by $\sigma_i(x)=\min(x_i,d-1)$ for all $x\in\Q{\ell}^n$ and $i\in [n]$; $\sigma$ is a surjection since $\ell\geq d$, and $\sigma$ acts as the identity on $\Q{d}^n$. Let $f\in F(n,\ell d)$ defined by 
\[
f(x,y)=(0^n,\sigma(x)). 
\] 
$f$ is 2-nilpotent because $f^2(x,y) = f(0^n,\sigma(x)) = (0^n,0^n)$. We will now prove that $\G(f)$ contains all the $d$-permutable digraphs on $[n]$. So let $G$ be a $d$-permutable digraph on $[n]$ and let $g\in F(G,d)$ be a permutation. We define $h\in F(n,\ell d)$ by 
\[
h(x,y)=(0^n,g(\sigma(x))). 
\] 
$h$ is 2-nilpotent because $h^2(x,y) = h(0^n,g(\sigma(x))) = (0^n,g(0^n))$. Furthermore, setting 
\[
\Omega=\{0^n\}\times\Q{d}^n,
\]
we have $\Im(f)=\Omega$ since $\sigma$ is a surjection, and $\Im(h)=\Omega$ since $g$ is a permutation of $\Q{d}^n$. Let $h'$ be the restriction of $h$ on  $\Q{d}^n\times\Q{d}^n$; since $\Omega\subseteq \Q{d}^n\times\Q{d}^n$ we have $h'\in F(n,d^2)$.  We now prove the following subgraphs relations:

\begin{itemize}
\item[(1)] $G(g)\subseteq G(h')$ and $G(h)\subseteq G(g)$.

\medskip
Suppose that $G(g)$ has an arc from $j$ to $i$. Then there exists $x,y\in\Q{d}^n$ which only differ in $x_j\neq y_j$ and such that $g_i(x)\neq g_i(y)$. Then $(x,x)$ and $(y,y)$ belongs to $\Q{d}^n\times\Q{d}^n$, and they only differ in component $j$, that is, $(x_j,x_j)\neq (y_j,y_j)$ but $(x_k,x_k)\neq (y_k,y_k)$ for all $k\neq j$. Furthermore, since $\sigma$ acts as the identity on $\Q{d}^n$, we have 
\[
h'_i(x,x)=(0,g_i(\sigma(x)))=(0,g_i(x))\neq (0,g_i(y))=(0,g_i(\sigma(y)))=h'_i(y,y),
\]
thus $G(h')$ has an arc from $j$ to $i$. This proves that $G(g)\subseteq G(h')$. Suppose now that $G(h)$ has an arc from $j$ to $i$. Let $(x,x'),(y,y')\in \Q{\ell }^n\times\Q{d}^n$ which only differ in component $j$ such that $h_i(x,x')\neq h_i(y,y')$. Since $h_i(x,x')=(0,g_i(\sigma(x)))$ and $h_i(y,y')=(0,g_i(\sigma(y)))$ we have $g_i(\sigma(x))\neq g_i(\sigma(y))$, thus $\sigma(x)\neq \sigma(y)$, and since $x$ and $y$ only differ in $x_j\neq y_j$ we deduce that $\sigma(x)$ and $\sigma(y)$ only differ in $\sigma(x)_j\neq \sigma(y)_j$. Thus $G(g)$ has an arc from $j$ to~$i$. This proves that $G(h)\subseteq G(g)$. 
\end{itemize}

\medskip
Since $h'$ is a restriction of $h$, we have $G(h')\subseteq G(h)$ and we deduce from (1) that $G(h)=G(g)$. It remains to prove that $h\sim f$. Let $\pi\in F(n,\ell d)$ defined by 
\[
\pi(x,y)=
\left\{
\begin{array}{ll}
(0^n,g(y))&\textrm{if $(x,y)\in\Omega$},\\
(x,y)&\textrm{otherwise}.
\end{array}
\right.
\]
Since $g$ is a permutation we have $\pi(\Omega)=\Omega$ and thus $\pi$ is a permutation. For $(x,y)\in\Omega$, we have 
\[
\pi(f(x,y))=\pi(0^n,0^n)=(0^n,g(0^n))=h(0^n,g(y))=h(\pi(x,y)). 
\]
For $(x,y)\not\in \Omega$, we have 
\[
\pi(f(x,y))=\pi(0^n,\sigma(x))=(0^n,g(\sigma(x)))=h(x,y)=h(\pi(x,y)). 
\]
Thus $\pi$ is an isomorphism between $h$ and $f$. 

\medskip
We now assume that $d<\ell$ and prove that $F(n,\ell d-1)$ also contains a $2$-nilpotent function which can be produced by all the $d$-permutable digraphs on $[n]$. Let $\Lambda$ be the set of $(x,y)\in \Q{\ell }^n\times\Q{d}^n$ such that $(x_i,y_i)\neq (\ell -1,d-1)$ for all $i\in [n]$. Thus $\Lambda$ is a subset of $\Q{\ell }^n\times\Q{d}^n$ of size $(\ell d-1)^n$ which can be identified with $\Q{\ell d-1}^n$. Let $f'$ and $h''$ be the restriction of $f$ and $h$ on $\Lambda$. Since $\Omega\subseteq \Lambda$, we have $f',h''\in F(n,\ell d-1)$. 

\medskip
Let us prove that $\G(f')$ contains all the $d$-permutable digraphs on $[n]$. For that it is sufficient to prove that $G(h'')=G(g)$ and $h''\sim f'$. Since $\ell >d$ we have $\Q{d}^n\times \Q{d}^n\subseteq \Lambda$. Hence $h'$ is a restriction of $h''$, and since $h''$ is a restriction of $h$, we have 
\[
G(h')\subseteq G(h'')\subseteq G(h).
\]
We then deduce from (1) that $G(h'')=G(g)$. It remains to prove that $h''\sim f'$. Let $\pi'$ be the restriction on $\pi$ on $\Lambda$. Since $\Omega\subseteq \Lambda$, $\pi'$ is a permutation of $\Lambda$. We then prove as above that $\pi'$ is an isomorphism between $h''$ and $f'$.  
\end{proof}

We now prove Lemma~\ref{lem:2-permutable}, that almost all digraphs are $2$-permutable. We first give a sufficient condition for $G$ to be $2$-permutable, and the proof of Lemma~\ref{lem:2-permutable} then consists in showing, by a simple counting argument,  that almost all digraphs satisfy this condition. Gadouleau \cite{G18rank} proved that the complete digraphs $K_n$ (with $n(n-1)$ arcs and no loop) is $2$-permutable, and Lemma~\ref{lem:2-permutable} is obtained by pushing forward his arguments.

\begin{lemma}\label{lem:2-permutable_condition}
Let $G$ be a coverable digraph on $[n]$ and, for each $i\in [n]$, let $N_i$ be the set of in-neighbors of $i$ in $G$. Suppose that, for all $1\leq i<j\leq n$, either $N_i\cap N_j=\emptyset$ or $|N_i\cap N_j|\geq (\log j)+2$. Then $G$ is $2$-permutable.
\end{lemma}

\begin{proof}
Since $G$ is coverable, there exists a permutation $\pi$ of $[n]$ such that $G$ has an arc from $\pi(i)$ to $i$ for all $i\in [n]$. For all $i,j\in [n]$, we set $M_i= N_i\setminus\{\pi(i)\}$ and $M_{ij}=M_i\cap M_j$. By hypothesis, if $i<j$, either $N_i\cap N_j=\emptyset$ or $|M_{ij}|\geq |N_i\cap N_j|-2\geq \log j$.  

\medskip
We first prove, by induction on $k=1,\dots,n$, that we can assign to each $i\in [k]$ a configuration $\alpha(i)\in\B^{M_i}$ in such a way that $\alpha(i)_{M_{ij}}\neq\alpha(j)_{M_{ij}}$ for every $1\leq i<j\leq k$ such that $M_{ij}\neq\emptyset$. For $k=1$ there is nothing to prove. So suppose that the property holds for $k<n$ and, to complete the induction, let us prove that, for $j=k+1$, there exists $\alpha(j)\in\B^{M_j}$ such that $\alpha(i)_{M_{ij}}\neq\alpha(j)_{M_{ij}}$ for all $1\leq i<j$ with $M_{ij}\neq \emptyset$. Let $I$ be the set of $1\leq i<j$ such that $M_{ij}\neq\emptyset$. For all $i\in I$, let $X_i$ be the set of configurations $\alpha\in\B^{M_j}$ which extends $\alpha(i)_{M_{ij}}$, that is, with $\alpha_{M_{ij}}=\alpha(i)_{M_{ij}}$. We have 
\[
|X_i|= 2^{|M_j|-|M_{ij}|}\leq 2^{|M_j|-\log j}=\frac{1}{j}2^{|M_j|}.
\]
Setting $X=\bigcup_{i\in I} X_i$ we obtain 
\[
|X|\leq \sum_{i\in I}|X_i|\leq \frac{|I|}{j}2^{|M_j|}\leq  \frac{j-1}{j}2^{|M_j|}<2^{|M_j|}.
\]  
Thus there exists $\alpha(j)\in \B^{|M_j|}\setminus X$, and then $\alpha(j)_{M_{ij}}\neq \alpha(i)_{M_{ij}}$ for all $i\in I$, as desired. 

\medskip
So let us assign to each $i\in [n]$ a configuration $\alpha(i)\in\B^{M_i}$ such that $\alpha(i)_{M_{ij}}\neq\alpha(j)_{M_{ij}}$ for every $1\leq i<j\leq n$ such that $M_{ij}\neq\emptyset$. Let $f\in F(n,2)$ defined by: for all $x\in\B^n$ and $i\in [n]$,  
\[
f_i(x)=\left\{
\begin{array}{ll}
x_{\pi(i)}&\textrm{if } x_{M_i}=\alpha(i),\\
1-x_{\pi(i)}&\textrm{otherwise}.
\end{array}
\right.
\]
It is clear that $f_i$ depends on $j$ if and only if $j\in M_i\cup\{\pi(i)\}=N_i$, and thus $G(f)=G$. Hence, it is sufficient to prove that $f$ is a permutation. So let $x,y\in\B^n$ distinct, and let us prove that $f(x)\neq f(y)$. Suppose, for a contradiction, that $f(x)=f(y)$. Since $x\neq y$, there exists $i\in [n]$ such that $x_{\pi(i)}\neq y_{\pi(i)}$; we take $i$ minimal for that property. Since $f_i(x)=f_i(y)$ we have $x_{M_i}=\alpha(i)\neq y_{M_i}$ or $x_{M_i}\neq \alpha(i)= y_{M_i}$. Suppose, without loss, that 
\[
x_{M_i}=\alpha(i)\neq y_{M_i}.
\]
Then there exists $j\in [n]$ such that $\pi(j)\in M_i$ and $x_{\pi(j)}\neq y_{\pi(j)}$; take $j$ minimal for that property. Since $\pi(i)\not\in M_i$ we have $i\neq j$ and, by the choice of $i$ we have $i<j$. Furthermore, since $\pi(j)\in N_i\cap N_j\neq\emptyset$ we have $|M_{ij}|\geq \log j\geq 1$ by hypothesis. Since $f_j(x)=f_j(y)$ we have $x_{M_j}=\alpha(j)\neq y_{M_j}$ or $x_{M_j}\neq \alpha(j)= y_{M_j}$. If $x_{M_j}=\alpha(j)$ then $\alpha(j)_{M_{ij}}=x_{M_{ij}}=\alpha(i)_{M_{ij}}$ and we obtain a contradiction since $M_{ij}\neq\emptyset$. Thus 
\[
x_{M_j}\neq \alpha(j)= y_{M_j}.
\]
We deduce that $x_{M_{ij}}\neq y_{M_{ij}}$ since otherwise $\alpha(i)_{M_{ij}}=x_{M_{ij}}=y_{M_{ij}}=\alpha(j)_{M_{ij}}$, and we obtain the same contradiction. Thus there exists $k\in [n]$ such that $\pi(k)\in M_{ij}$ and $x_{\pi(k)}\neq y_{\pi(k)}$. Since $\pi(i),\pi(j)\not\in M_{ij}$, we have $k\neq i,j$ and since $\pi(k)\in M_i$, by the choice of $j$, we have $j<k$. Furthermore, since $\pi(k)\in M_{ij}$, we have $N_i\cap N_k\neq\emptyset$ and $N_j\cap N_k\neq\emptyset$ thus $|M_{ik}|,|M_{jk}|\geq \log k\geq 1$ by hypothesis. Since $f_k(x)=f_k(y)$ we have $x_{M_k}=\alpha(k)\neq y_{M_k}$ or $x_{M_k}\neq \alpha(k)= y_{M_k}$. If $x_{M_k}=\alpha(k)$ then $\alpha(k)_{M_{ik}}=x_{M_{ik}}=\alpha(i)_{M_{ik}}$ and we obtain a contradiction since $M_{ik}\neq\emptyset$. If $y_{M_k}=\alpha(k)$ then $\alpha(k)_{M_{jk}}=y_{M_{jk}}=\alpha(j)_{M_{jk}}$ and we obtain a contradiction since $M_{jk}\neq\emptyset$. Thus $f$ is indeed a permutation. 
\end{proof}

\begin{proof}[\BF{Proof of Lemma~\ref{lem:2-permutable}}]
Let $\mathcal{G}(n)$ be the set of digraphs $G$ on $[n]$ such that, for all distinct $i,j\in [n]$, $i$ and $j$ have at least $\log n+2$ common in-neighbors. By Lemma~\ref{lem:2-permutable_condition}, any coverable digraph in $\mathcal{G}(n)$ is $2$-permutable. Since almost all digraphs are Hamiltonian, and thus coverable, it is sufficient to prove that $|\mathcal{G}(n)|/2^{n^2}\to 1$ as $n\to\infty$.

\medskip
Let $\lambda(n)=2^{n^2}-|\mathcal{G}(n)|$. Denoting $\TWO^{[n]}$ the set of subsets of $[n]$, and setting $\ell(n)=\lfloor \log n+1\rfloor$, we see that $\lambda(n)$ is exactly the number of tuples $(N_1,\dots,N_n)\in (\TWO^{[n]})^n$ such that $|N_i\cap N_j|\leq\ell(n)$ for some distinct $i,j\in [n]$. Let $\lambda_2(n)$ be the number of couples $(N,N')\in (\TWO^{[n]})^2$ such that $|N\cap N'|\leq\ell(n)$. For $N\in \TWO^{[n]}$ fixed of size $k$, there are $\sum_{i=0}^{\ell(n)}{k\choose i}2^{n-k}$ sets $N'\in \TWO^{[n]}$ such that $|N\cap N'|\leq\ell(n)$. We deduce that 
\begin{eqnarray*}
\lambda_2(n)&=&\sum_{k=0}^n {n\choose k}\sum_{i=0}^{\ell(n)}{k\choose i}2^{n-k}\\
&\leq& \sum_{k=0}^n {n\choose k}(k+1)^{\ell(n)} 2^{n-k}\\
&\leq& (n+1)^{\ell(n)} \sum_{k=0}^n {n\choose k}2^{n-k}\\
&=& (n+1)^{\ell(n)} 3^n.
\end{eqnarray*}
For $i,j\in [n]$ distinct, let $\lambda_{ij}(n)$ be the number of tuples $(N_1,\dots,N_n)\in (\TWO^{[n]})^n$ such that $|N_i\cap N_j|\leq\ell(n)$. We have $\lambda_{ij}(n)=\lambda_2(n)(2^n)^{n-2}\leq (n+1)^{\ell(n)} 3^n(2^n)^{n-2}$ and thus, setting  $c=2-\log 3$, we obtain  
\[
\lambda(n)\leq \sum_{i,j\in [n]\atop i\neq j} \lambda_{ij}(n)\leq n^2(n+1)^{\ell(n)} 3^n 2^{n(n-2)}\leq 2^{n^2-cn+2\ell(n)\log 2n}\leq 2^{n^2-cn+2(\log 2n)^2}.
\] 
Thus
\[
\frac{|\mathcal{G}(n)|}{2^{n^2}}=1-\frac{\lambda(n)}{2^{n^2}}\geq 1-2^{-cn+2(\log 2n)^2}.
\]
Since $c>0$ we have $2^{-cn+2(\log 2n)^2}\to 0$ as $n\to\infty$. Thus $|\mathcal{G}(n)|/2^{n^2}\to 1$ as $n\to\infty$. 
\end{proof}

\section{Binary case}\label{sec:binary}

We were not able to give any non trivial information on the asymptotic behavior of $\gamma_2(n)$, and that $\gamma_2(n)/2^{n^2}$ tends to $0$ or $1$ is open. Instead, we present some weak bounds on $\gamma_2(n)$, which however highlight some specificities of the binary case, and then present a modified version of universality which is more easy to analyze. 

\medskip
We begin by showing that Theorem \ref{thm:Hamiltonian2} fails in the binary case.  

\begin{proposition}
For all $n\geq 2$, there does not exist $f\in F(n,2)$ such that $\G(f)$ contains all the Hamiltonian digraphs on $[n]$. 
\end{proposition}

\begin{proof}
Let $f\in F(n,2)$, and let $G$ be any Hamiltonian digraph on $[n]$ where each vertex is of in-degree $2$. We will prove that $\G(f)$ cannot contain both $G$ and $C_n$. Suppose for a contradiction, that $\G(f)$ contains both $C_n$ and $G$. Since $F(C_n,2)$ only contains permutations (Lemma \ref{lem:disjoint_union_of_cycles}), $f$ is a permutation, and to obtain a contradiction, it is sufficient to prove that $F(G,2)$ does not contain any permutation. Suppose, for a contradiction that $h\in F(G,2)$ is a permutation. Then, for all $i\in [n]$, we have $|h^{-1}_i(0)|=|h^{-1}_i(1)|=2^{n-1}$. Since $i$ is of in-degree $2$, this forces $h_i$ to be affine, that is, there is $a_i\in\B$ such that $h_i(x)=x_{i_1}\oplus x_{i_2}\oplus a_i$ for all $x\in\B^n$, where $i_1$ and $i_2$ are the two in-neighbors of $i$. But then, $h(0^n)=h(1^n)=(a_1,\dots,a_n)$, a contradiction. 
\end{proof}

\subsection{Lower bounds}

We now exhibit our best construction of a function $f\in F(n,2)$ with $|\G(f)|$ large. Given a digraph $H$ on $[n-1]$, the \EM{universal augmentation} of $H$ is the digraph $G$ on $[n]$ obtained from $H$ by adding vertex $n$ and an arc from $n$ to every $i\in [n-1]$. 

\begin{proposition}\label{pro:universal_augmentation}
For all $n\geq 2$, there exists a $2$-nilpotent function $f\in F(n,2)$ such that $\G(f)$ contains the universal augmentation of the all the $2$-permutable digraphs on $[n-1]$. 
\end{proposition}

\begin{proof}
Let $f\in F(n,2)$ defined as follows: for all $x\in\B^{n-1}$, 
\[
f(x,1)=(x,0)\quad\textrm{and}\quad f(x,0)=(0^{n-1},0).
\]
We easily check that $f$ is $2$-nilpotent, with fixed point $0^n$. We will prove that $\G(f)$ has the desired property. So let $H$ be a $2$-permutable digraph on $[n-1]$, and let $g\in F(H,2)$ be a permutation. Let $h\in F(n,2)$ defined as follows: for all $x\in\B^{n-1}$, 
\[
h(x,1)=(g(x),0)\quad\textrm{and}\quad h(x,0)=(g(0^{n-1}),0).
\]
We easily check that $H$ is obtained from $G(h)$ by deleting vertex $n$. Furthermore, $G(h)$ has an arc from $n$ to each $i\in [n-1]$: since $g$ is a permutation, $g_i\neq\cst$, thus there exists $x\in\B^{n-1}$ such that $g_i(x)\neq g_i(0^{n-1})$ and thus $h_i(x,1)=g_i(x)\neq g_i(0^{n-1})=h_i(x,0)$. Since $h_n=\cst=0$, $n$ is a source of $G(h)$ and thus $G(h)$ is the universal augmentation of $H$. It remains to prove that $h\sim f$.  Let $\pi\in F(n,2)$ defined as follows:  for all $x\in\B^{n-1}$, 
\[
\pi(x,1)=(x,1)\quad\textrm{and}\quad \pi(x,0)=(g(x),0).
\]
Then $\pi$ is a permutation of $\B^n$ since $g$ is a permutation. Furthermore, 
\[
\pi(f(x,1))=\pi(x,0)=(g(x),0)=h(x,1)=h(\pi(x,1))
\]
and
\[
\pi(f(x,0))=\pi(0^{n-1},0)=(g(0^{n-1}),0)=h(g(x),0)=h(\pi(x,0)).
\]
Thus $h$ is indeed isomorphic to $f$. 
\end{proof}

The corresponding lower bound on $\gamma_2(n)$ is the following. 

\begin{proposition}
For all $\epsilon>0$ and $n$ large enough, $\gamma_2(n)\geq (1-\epsilon)2^{(n-1)^2}$. 
\end{proposition}

\begin{proof}
Let $p(n)$ be the number of $2$-permutable digraphs on $[n]$. By Proposition \ref{pro:universal_augmentation} we have $\gamma_2(n)\geq p(n-1)$. By Lemma~\ref{lem:2-permutable}, almost all digraphs are $2$-permutable. This is equivalent to say that, for all $\epsilon>0$, there exists $n_0$ such that $p(n)\geq (1-\epsilon)2^{n^2}$ for all $n\geq n_0$. Consequently, $\gamma_2(n)\geq (1-\epsilon)2^{(n-1)^2}$ for all $n>n_0$. 
\end{proof}

\subsection{Upper bounds}

To get an upper bound on $\gamma_2(n)$, we exhibit two families of digraphs that cannot cohabit in $\G(f)$. Let us say that a digraph $G$ has an \EM{initial loop} if it has a vertex $i$ such that $i$ has a unique in-neighbor, which is itself. 

\begin{proposition}\label{pro:initial_loop}
For all $n\geq 1$, there does not exist $f\in F(n,2)$ such that $\G(f)$ contains both an acyclic digraph and a digraph with an initial loop. 
\end{proposition}

\begin{proof}
Suppose that $\G(f)$ contains an acyclic digraph. Then $f^n=\cst$ by Robert's Theorem (Theorem \ref{thm:robert}), and thus $f^{2n}=\cst$. Suppose now, for a contradiction, that $\G(f)$ contains a digraph $G$ with an initial loop, say on vertex $i$, and let $h\in F(G,2)$ isomorphic to $f$. Since $i$ has a unique in-neighbor, which is itself, there is $a\in\B$ such that $h_1(x)=x_1\oplus a$ for all $x\in\B^n$. Consequently, for every integer $p\geq 0$ and $x\in\B^n$, we have $h^{2p}_1(x)=x_1$. We deduce that $h^{2n}\neq \cst$, and thus $f^{2n}\neq\cst$, a contradiction. 
\end{proof}

The corresponding upper bound on $\gamma_2(n)$ is the following. 

\begin{proposition}
For all $n\geq 2$, $\gamma_2(n)\leq 2^{n^2}- 2^{\frac{n^2}{2}}$.
\end{proposition}

\begin{proof}
Let $a(n)$ be the number of acyclic digraphs on $[n]$, and let $b(n)$ be the number of digraphs on $[n]$ with an initial loop. By Proposition \ref{pro:initial_loop} we have 
\[
|\G(f)|\leq 2^{n^2}-\min(a(n),b(n)).
\]

\medskip
Given $i\in [n]$, the number of digraphs on $[n]$ with an initial loop on $i$ is exactly $(2^n)^{n-1}$, and thus $b(n)\geq (2^n)^{n-1}$. Let us now estimate $a(n)$. Let $\pi$ be a permutation of $[n]$. Let $A(n,\pi)$ be the set of acyclic digraphs on $[n]$ with the Hamiltonian path $\pi(1),\pi(2),\dots,\pi(n)$. If $\sigma$ is a permutation of $[n]$ distinct from $\pi$ then $A(n,\pi)\cap A(n,\sigma)=\emptyset$, thus $a(n)\geq n!|A(n,\pi)|$. Let $T$ be the transitive tournament contained in $A(n,\pi)$. It has ${n\choose 2}-(n-1)$ arcs not contained in the Hamiltonian path. Any digraph obtained from $T$ by deleting some of these arcs is in $A(n,\pi)$ and thus 
\[
A(n,\pi)\geq 2^{{n\choose 2}-n+1}.
\]
We obtain  
\[
a(n)\geq n!2^{{n\choose 2}-n+1}\geq 2^{\frac{n^2}{2}}.
\]
Since $b(n)\geq (2^n)^{n-1}\geq 2^{\frac{n^2}{2}}$ we obtain the bound of the statement. 
\end{proof}

\subsection{Induced universality}

We finally introduce a modified notion of universality, which is more easy to analyze. Given an integer $k\geq 1$, we say that $f\in F(n,2)$ is \EM{$k$-induced universal} if, for every $k$-vertex digraph $H$, some members of $\G(f)$ have an induced subgraph isomorphic to $H$. Since $\G(f)$ is closed under isomorphism, this is equivalent to say that, for every digraph $H$ on $[k]$, there exists $G\in\G(f)$ such that $H$ is the subgraph of $G$ induced by $[k]$. For fixed $k\geq 1$, we are interested in the smallest integer \EM{$n(k)$} such that $F(n(k),2)$ contains a $k$-induced universal function. We first prove that $n(k)$ is  a threshold for the presence of a $k$-induced universal function. 

\begin{lemma}
$F(n,2)$ contains a $k$-induced universal function iff $n\geq n(k)$. 
\end{lemma}

\begin{proof}
By definition, if $n<n(k)$ then $F(n,2)$ does not contain a $k$-induced universal function. For the converse, we assume that $F(n,2)$ contains a $k$-induced universal function $f$, and we prove that $F(n+1,2)$ also contains such a function. Indeed, let $f'\in F(n+1,2)$ defined by: for all $x\in\B^n$ and $a\in\B$, $f'(x,a)=(f(x),a)$. Let $H$ be a digraph on $[k]$, and let $h\in F(n,2)$ isomorphic to $f$ such that $H$ is an induced subgraph of $G(h)$. Let $h'\in F(n+1,2)$ defined by: for all $x\in\B^n$ and $a\in\B$, $h'(x,a)=(h(x),a)$. Since $h\sim f$ we have $h'\sim f'$. Furthermore, $G(h')$ is obtained from $G(h)$ by adding the vertex $n+1$ and a loop on it. Thus $H$ is induced subgraph of $G(h')$ and we deduce that $f'$ is $k$-induced universal.
\end{proof}

We now estimate $n(k)$ rather accurately. 

\begin{theorem}\label{thm:Nk}
For all $k\geq 1$, 
\[
k+\log k-\log(\lceil \log k\rceil+1)-1\leq n(k)\leq k+\lceil \log k\rceil+1. 
\]
\end{theorem}  

We start with the upper bound. 

\begin{lemma}\label{lem:Nk_upper_bound}
For all $k\geq 1$, we have $n(k)\leq k+\lceil \log k\rceil+1$. 
\end{lemma}

\begin{proof}
Let $n=k+\lceil \log k\rceil+1$ and let $m=n-k-1$. We identify $\B^n$ with $\B^{k+1}\times \B^m$. Let $f\in F(n,2)$ defined by 
\[
f(x,y)=
\left\{
\begin{array}{ll}
(0^{k+1},y)&\textrm{if $x_1=0$}\\
(1^{k+1},y)&\textrm{if $x_1=1$}.
\end{array}
\right.
\]
Note that $f$ has exactly $2^{m+1}$ fixed points, which are $(0^{k+1},y)$ and $(1^{k+1},y)$ for each $y\in\B^m$, and each fixed point has exactly $2^k$ pre-images. We will prove that $f$ is $k$-induced universal. So let $H$ be a digraph on $[k]$, and let us prove that some member of $\G(f)$ has an induced subgraph isomorphic to $H$. 

\medskip
Since $m=\lceil \log k\rceil$, there exists a surjection $\phi$ from $\B^m$ to $[k]$. For each $i\in [k]$, let $N_i$ be the out-neighbors of $i$ in $H$. We then define $I_i,J_i$ as follows: if $i\not\in N_i$, then $I_i=\{i,k+1\}$ and $J_i=N_i\cup\{k+1\}$; and if $i\in N_i$, then $I_i=\{i\}$ and $J_i=N_i$. Note that, in both cases, $|I_i\cap J_i|=1$. Let $h\in F(n,2)$ defined by 
\[
h(x,y)=
\left\{
\begin{array}{ll}
(0^{k+1},y)&\textrm{if $\sum_{i\in I_{\phi(y)}}x_i=0\mod 2$}\\
(e_{J_{\phi(y)}},y)&\textrm{if $\sum_{i\in I_{\phi(y)}}x_i=1\mod 2$}.
\end{array}
\right.
\]
 
\medskip
For every $y\in\B^m$, it is clear that $(0^n,y)$ is a fixed point of $h$, and since $|I_{\phi(y)}\cap J_{\phi(y)}|$ is odd, $(e_{J_{\phi(y)}},y)$ is also a fixed point. Hence $h$ has exactly $2^m$ fixed points, and we easily check that each of them has exactly $2^k$ pre-images. We deduce that $h\sim f$, so to prove that some member of $\G(f)$ has an induced subgraph isomorphic to $H$, it is sufficient to prove that $H$ is the subgraph of $G(h)$ induced by $[k]$. 

\medskip
Suppose that $H$ has an arc from $i$ to $j$. Since $\phi$ is a surjection, there exists $y\in\B^m$ with $\phi(y)=i$. Then $h(0^{k+1},y)=(0^{k+1},y)$ and, since $i\in I_i$, we have $h(e_i,y)=(e_{J_i},y)$. Since $j\in J_i$ we deduce that $G(h)$ has an arc from $i$ to $j$. So $H$ is a subgraph of $G(h)$. To prove that it is an induced subgraph, suppose that $G(h)$ has an arc from $i$ to $j$ with $i,j\in [k]$. Then there exists two configurations $x,x'\in\B^{k+1}$, which only differ in component $i$, and a configuration $y\in\B^m$, such that $h(x,y)$ and $h(x',y)$ differ in component $j$. This obviously implies $j\in J_{\phi(y)}$ and $i\in I_{\phi(y)}$, and thus $j\in N_i$. Consequently, $H$ is the subgraph of $G(h)$ induced by $[k]$. 
\end{proof}

For the lower bound we will use the following two lemmas concerning the rank of a function $f\in F(n,2)$. The first says that a large independent set in $G(f)$ forces a small rank, generalizing the trivial observation that if $G(f)$ has no arc (independent set of size $n$) then $f$ has rank one. The second says that a large induced subgraph in $G(f)$ containing a loop on each vertex and no other arcs forces a large rank, generalizing the trivial observation that if $G(f)=L_{n,n}$ then $f$ has rank $2^n$.

\begin{lemma}\label{lem:independent_set}
Let $f\in F(n,2)$ and suppose that $G(f)$ has an independent set of size $k$. Then the rank of $f$ is at most $2^{2(n-k)}$. 
\end{lemma}

\begin{proof}
Let $I$ be an independent set of $G(f)$ of size $k$, and $J=[n]\setminus I$. For any $x\in\B^n$, the restriction $f(x)_I$ only depends on $x_J$, and thus there is at most $2^{|J|}$ such restrictions. Since they are obviously at most $2^{|J|}$ restrictions $f(x)_J$, and since $(I,J)$ is a partition of $[n]$, we deduce that the number of images is at most $2^{|J|}\times 2^{|J|}=2^{2(n-k)}$.
\end{proof}

\begin{lemma}\label{lem:kC1}
Let $f\in F(n,2)$ and suppose that $G(f)$ has an induced subgraph isomorphic to $L_{k,k}$. Then the rank of $f$ is at least $2^{k/2^{n-k}}$. 
\end{lemma}

\begin{proof}
Let $I\subseteq [n]$ such that the subgraph of $G(f)$ induced by $I$ is isomorphic to $L_{k,k}$. It means that, for each $i\in I$, $G(f)$ has no arc from $I\setminus\{i\}$ to $i$, and it has a loop on $i$. So there is a configuration $x^i$ such that $f_i(x^i)\neq f_i(x^i+e_i)$, and without loss, $x^i_i=0$. Let $J=[n]\setminus I$. Note that $|I|=k$ and $|J|=n-k$. Considering the function $i\to x^i_J$, by the pigeonhole principle, there exists $z\in\B^J$ and a subset $K\subseteq I$ of size at least $k/2^{n-k}$ such that $x^i_J=z$ for all $i\in K$. Let $X$ be the set of $x\in\B^n$ with $x_J=z$ and $x_i=0$ for all $i\in I\setminus K$. Then $|X|=2^{|K|}\geq 2^{k/2^{n-k}}$. Thus it is sufficient to prove that $f$ is injective on $X$. So let distinct $x,y\in X$. Then there exists $i\in K$ such that $x_i\neq y_i$, say $x_i<y_i$ without loss. Since $x^i_i=x_i=0$ and $x^i_J=z=x_J$, and since $G(f)$ has no arc from $I\setminus \{i\}$ to $i$, we have $f_i(x^i)=f_i(x)$. Similarly, $(x^i+e_i)_i=y_i=1$ and $(x^i+e_i)_J=z=y_J$ and thus $f_i(x^i+e_i)=f_i(y)$. We deduce that $f_i(x)\neq f_i(y)$ and thus $f(x)\neq f(y)$ as desired. 
\end{proof}

We are now in a position to prove the lower bound in Theorem~\ref{thm:Nk}.

\begin{lemma}\label{lem:Nk_lower_bound}
For all $k\geq 1$, we have $n(k)\geq k+\log k-\log(\lceil \log k\rceil+1)-1$.
\end{lemma}

\begin{proof}
Let $n=n(k)$, $m=n(k)-k$, and let $f\in F(n(k),2)$ be a $k$-induced universal function. Then there exists $g,h\in F(n,2)$ isomorphic to $f$ such that $G(g)$ has an independent set of size $k$ and $G(h)$ has an induced subgraph isomorphic to $L_{k,k}$. Since $f,g,h$ have the same rank, we deduce from Lemmas~\ref{lem:independent_set} and \ref{lem:kC1} that the rank of $f$ is at most $2^{2m}$ and at least $2^{k/2^m}$. So $2^{2m}\geq 2^{k/2^m}$ that is, $m\geq \log k-\log m-1$ and thus $n(k)\geq k+\log k-\log m-1$. By Lemma~\ref{lem:Nk_upper_bound} we have $n(k)\leq k+\lceil \log k\rceil+1$. So $m\leq \lceil \log k\rceil+1$ and thus $n(k)\geq k+\log k-\log(\lceil \log k\rceil+1)-1$.
\end{proof}

\section{Concluding remarks}\label{sec:conclu}

\begin{itemize}
\item 
The main open problem raised in this paper is whether there exists almost universal functions in the binary case: do we have $\gamma_2(n)/2^{n^2}\to 1$ as $n\to\infty$? 
\item
Is it difficult to decide if a given function $f\in F(n,q)$ is universal? By Theorem \ref{thm:universal_1}, it is sufficient to check that $f$ is $2$-nilpotent and that $L_{1,n},L_{n,n}\in\G(f)$. We can obviously check in polynomial time if $f$ is $2$-nilpotent and, by Lemma~\ref{lem:L1n}, we can also check if $L_{1,n}\in\G(f)$ in polynomial time (the size of the input $f$ being $q^n$). So what is the complexity of deciding if $L_{n,n}\in\G(f)$? More generally, given a fixed digraph $H$, what is the complexity of deciding if $H\in\G(f)$? This problem is in NP, but is it NP-complete for some $H$?  
\item 
Generalizing the previous point, given a fixed family $\mathcal{H}$ of digraphs, what is the complexity of deciding if $\G(f)$ intersects $\mathcal{H}$? This seems particularly interesting when $\mathcal{H}$ is the family of acyclic digraphs $\mathcal{A}$.  By Robert's theorem, if $\G(f)$ intersects $\mathcal{A}$ then $f$ has a simple structure: $f$ is $k$-nilpotent for some $k\leq n$. But is there some hidden complexity: is it NP-complete to decide if $\G(f)$ intersects $\mathcal{A}$?  

\end{itemize}

\paragraph{Acknowledgments} This work has been partially funded by the HORIZON-MSCA-2022-
SE-01 project 101131549 ‘‘ACANCOS’’ project and the
ANR-24-CE48-7504 ‘‘ALARICE’’ project.


\bibliographystyle{plain}
\bibliography{BIB}

\end{document}